\documentclass[12pt]{amsart}
\usepackage{amsbsy}
\usepackage{graphicx,epsfig,subfigure,psfrag}
\begin{document}

\newtheorem{theorem}{Theorem}
\newtheorem{proposition}{Proposition}
\newtheorem{lemma}{Lemma}
\newtheorem{corollary}{Corollary}
\newtheorem{definition}{Definition}
\newtheorem{remark}{Remark}
\newcommand{\tex}{\textstyle}
\numberwithin{equation}{section} \numberwithin{theorem}{section}
\numberwithin{proposition}{section} \numberwithin{lemma}{section}
\numberwithin{corollary}{section}
\numberwithin{definition}{section} \numberwithin{remark}{section}
\newcommand{\ren}{\mathbb{R}^N}
\newcommand{\re}{\mathbb{R}}
\newcommand{\n}{\nabla}
\newcommand{\iy}{\infty}
\newcommand{\pa}{\partial}
\newcommand{\fp}{\noindent}
\newcommand{\ms}{\medskip\vskip-.1cm}
\newcommand{\mpb}{\medskip}
\newcommand{\AAA}{{\bf A}}
\newcommand{\BB}{{\bf B}}
\newcommand{\CC}{{\bf C}}
\newcommand{\DD}{{\bf D}}
\newcommand{\EE}{{\bf E}}
\newcommand{\FF}{{\bf F}}
\newcommand{\GG}{{\bf G}}
\newcommand{\oo}{{\mathbf \omega}}
\newcommand{\Am}{{\bf A}_{2m}}
\newcommand{\CCC}{{\mathbf  C}}
\newcommand{\II}{{\mathrm{Im}}\,}
\newcommand{\RR}{{\mathrm{Re}}\,}
\newcommand{\eee}{{\mathrm  e}}
\newcommand{\LL}{L^2_\rho(\ren)}
\newcommand{\LLL}{L^2_{\rho^*}(\ren)}
\renewcommand{\a}{\alpha}
\renewcommand{\b}{\beta}
\newcommand{\g}{\gamma}
\newcommand{\G}{\Gamma}
\renewcommand{\d}{\delta}
\newcommand{\D}{\Delta}
\newcommand{\e}{\varepsilon}
\newcommand{\var}{\varphi}
\renewcommand{\l}{\lambda}
\renewcommand{\o}{\omega}
\renewcommand{\O}{\Omega}
\newcommand{\s}{\sigma}
\renewcommand{\t}{\tau}
\renewcommand{\th}{\theta}
\newcommand{\z}{\zeta}
\newcommand{\wx}{\widetilde x}
\newcommand{\wt}{\widetilde t}
\newcommand{\noi}{\noindent}
\newcommand{\uu}{{\bf u}}
\newcommand{\xx}{{\bf x}}
\newcommand{\yy}{{\bf y}}
\newcommand{\zz}{{\bf z}}
\newcommand{\aaa}{{\bf a}}
\newcommand{\cc}{{\bf c}}
\newcommand{\jj}{{\bf j}}
\newcommand{\ggg}{{\bf g}}
\newcommand{\UU}{{\bf U}}
\newcommand{\YY}{{\bf Y}}
\newcommand{\HH}{{\bf H}}
\newcommand{\GGG}{{\bf G}}
\newcommand{\VV}{{\bf V}}
\newcommand{\ww}{{\bf w}}
\newcommand{\vv}{{\bf v}}
\newcommand{\hh}{{\bf h}}
\newcommand{\di}{{\rm div}\,}
\newcommand{\ii}{{\rm i}\,}
\newcommand{\inA}{\quad \mbox{in} \quad \ren \times \re_+}
\newcommand{\inB}{\quad \mbox{in} \quad}
\newcommand{\inC}{\quad \mbox{in} \quad \re \times \re_+}
\newcommand{\inD}{\quad \mbox{in} \quad \re}
\newcommand{\forA}{\quad \mbox{for} \quad}
\newcommand{\whereA}{,\quad \mbox{where} \quad}
\newcommand{\asA}{\quad \mbox{as} \quad}
\newcommand{\andA}{\quad \mbox{and} \quad}
\newcommand{\withA}{,\quad \mbox{with} \quad}
\newcommand{\orA}{,\quad \mbox{or} \quad}
\newcommand{\ef}{\eqref}
\newcommand{\ssk}{\smallskip}
\newcommand{\LongA}{\quad \Longrightarrow \quad}
\def\com#1{\fbox{\parbox{6in}{\texttt{#1}}}}
\def\N{{\mathbb N}}
\def\A{{\cal A}}
\newcommand{\de}{\,d}
\newcommand{\eps}{\varepsilon}
\newcommand{\be}{\begin{equation}}
\newcommand{\ee}{\end{equation}}
\newcommand{\spt}{{\mbox spt}}
\newcommand{\ind}{{\mbox ind}}
\newcommand{\supp}{{\mbox supp}}
\newcommand{\dip}{\displaystyle}
\newcommand{\prt}{\partial}
\renewcommand{\theequation}{\thesection.\arabic{equation}}
\renewcommand{\baselinestretch}{1.1}
\newcommand{\Dm}{(-\D)^m}

\title
{\bf  Blow-up in higher-order reaction-diffusion\\ and wave
equations: how
    $\sqrt{\mbox{log\,log}}$
 factor occurs}


\author {V.A.~Galaktionov}

\address{Department of Mathematical Sciences, University of Bath,
 Bath BA2 7AY, UK}
\email{vag@maths.bath.ac.uk}




\date{\today}

\begin{abstract}

The origin of non self-similar blow-up in higher-order
reaction-diffusion (parabolic) or wave (hyperbolic) equations with
typical models
 $$
  \tex{
 u_t= u^2(-u_{xxxx}+u) \quad
 \mbox{or} \quad  u_{tt}= u^2(-u_{xxxx}+u)\inB  (-L,L) \times (0,T),
 }
 $$
 with zero Dirichlet boundary conditions at $x= \pm L$, where $L>L_0\in( \frac \pi 2,
 \pi)$ is discussed. The rate of
 blow-up is shown to get an extra universal
 $\sqrt{\ln|\ln(T-t)|}$ factor in addition to the standard similarity
one $\frac 1{\sqrt{T-t}}$. The explanation is based on matching
with  the so-called {\em logarithmic travelling waves} as group
invariant solutions of the equation.

 Some links and similarities with double-log blow-up terms occurring in earlier studies of plasma
 physics parabolic equations and the nonlinear critical Schr\"odinger
equation  are discussed. On the other hand, $\sqrt{\ln|\ln(T-t)|}$
obtained in Petrovskii's boundary regularity study for the heat
equation in 1934 was the first its appearance in PDE theory.


\end{abstract}

\maketitle

\section{Introduction: {\bf THREE} mysterious blow-up
$\sqrt{\log\, \log}$ of PDE theory}
\label{S1}

\subsection{On our main goal}
In the middle of the 1980s, in the study of singularity formation
phenomena of blow-up in a reaction-diffusion equation from plasma
physics and, almost simultaneously, in self-focusing for the cubic
nonlinear Schr\"odinger equation, the physical and formal
asymptotic methods, to say nothing about rigorous justifications,
faced an extremely difficult issue of appearance the so-called
{\em double logarithmic}, or $\sqrt{\log\,\log}$ factor:
 \be
 \label{AA1}
 A_0(t) \sim \sqrt{\ln|\ln(T-t)|} \asA t \to T^-.
  \ee
 Here $T<\infty$ is the blow-up time in the sense that the
 solution $u=u(x,t)$ of the PDE under consideration is well-defined and is classic
for all $t \in (0,T)$, but gets unbounded\footnote{For a full
correctness, the sign ``$\limsup$" should be used; however, for
problems that are locally well-posed in $L^\iy$, evidently, just
``$\lim$" does.}:
 \be
 \label{AA2}
  \tex{
  \lim_{t \to T^-} \sup_x |u(x,t)| =+\iy.
  }
   \ee

The goal of this paper is to introduce a number of higher-order
scaling invariant nonlinear PDEs of parabolic, hyperbolic, and
nonlinear dispersion types, which can exhibit the factor \ef{AA1}
in their blow-up behaviour. In fact, we are going to address a
wider question:
 \be
 \label{Q1}
  \fbox{$
  \mbox{under which assumptions on  PDEs, (\ref{AA1}) occurs in blow-up or not,}
   $}
   \ee
 and how the absence of \ef{AA1} affects the generic blow-up
 behaviour.
The mysterious (becoming not that much after a proof is found)
blow-up factor \ef{AA1} is already well-known for a few PDEs, so
we are inevitably obliged to begin with this amazing history of
the twentieth century.

\subsection{FIRST $\sqrt{\bf log\,log}$ in classic
 parabolic theory: boundary regularity and Petrovskii's
criterion, 1934}

 It is truly amazing that
 its actual origin lies in the heart of
 PDE theory: {\em regularity of a boundary point}.
 For the Dirichlet problem for the Laplace equation, this study,
 began by Green, Gauss, Lord Kelvin, Dirichlet in the
 first half of the nineteenth century and by many other great mathematicians,
 was completed by Wiener in 1924  \cite{Wien24}, who derived his
 famous regularity criterion (a necessary and sufficient condition). A detailed history of potential theory can be found in
Kellogg \cite[pp.~277--285]{Kell29}.

The same  regularity question for the heat equation in a
non-cylindrical domain was in 1934 initiated by Petrovskii
\cite{Pet34Cr, Pet35}, where the double-log actually occurred for
first time.
This is the question on {\em irregular} or {\em regular} point
$(x,t)=(0,T)$ for the 1D heat equation
 \be
 \label{ir1}
  \left\{
 \begin{matrix}
   u_t=u_{xx} \quad \mbox{in}
 \quad Q_T=\{|x|<R(t), \,\,\, 0<t<T\} \whereA R(t) \to 0^+ \,\,\mbox{as}\,\, t \to
 T^-, \ssk\ssk\\
  \mbox{with bounded smooth data $u(x,0)=u_0(x)$ on $[-R(0),R(0)]$}.
  \qquad\qquad\qquad\quad\,\,\,\,\,
   \end{matrix}
    \right.
  \ee
Here the lateral boundary $\{x= \pm R(t)\}$ is given by a function
$R(t)$ that is  assumed to be positive and $C^1$-smooth for all $0
\le t <T$ and is allowed to have  a singularity of $R'(t)$ at
$t=T^-$ only. Then  the value of $u(x,t)$ is studied at the end
``blow-up" point $(0,T)$, to which the domain $Q_T$ ``shrinks" as
$t \to T^-$.
Thus, $(x,t)=(0,T)$ is {\em regular}, if any value of the solution
$u(x,t)$ can be prescribed there by continuity as a standard
boundary value on $\partial Q_T$. Otherwise the point is {\em
irregular}, if the value $u(0,T)$ is not arbitrary and is given by
the
evolution as $t \to T^-$.



 Petrovskii in 1934--35 \cite{Pet34Cr, Pet35}, in particular, established the following:
 \be
 \label{PP1}
 \fbox{$
 \begin{matrix}
{\rm (i)}\,R(t)=2 \sqrt{T-t} \,\,\sqrt{\ln|\ln (T-t)|} \LongA
\mbox{$(0,T)$ is regular}, \andA \qquad\quad \ssk\ssk\\ {\rm
(ii)}\, R(t)=2(1+\e)\sqrt{T-t}\,\, \sqrt{ \ln|\ln (T-t)|}, \,\,
\e>0 \LongA
 \mbox{$(0,T)$ is irregular};
  \end{matrix}
  $}
  \ee
see the most  recent survey in  \cite{GalPet2m}  for a full list
of references. As far as we know, \ef{PP1} this  is the first
clear appearance of the blow-up $\sqrt{\rm log\,log}$ factor in
PDE theory.
\subsection{SECOND $\sqrt{\bf log\,log}$ in the NLSE (1985): the origin and
 beginning of the
log-log story in nonlinear PDE theory}

We now return  to the true blow-up problems, but will begin with
another famous PDE, which is not a subject of the present paper,
but used to be well-recognized as the source of the log-log. It
has been well-accepted that  the origin (besides Petrovskii's
result more belonging to regularity and even probability theory)
of such an unusual rate of {\em blow-up} including the
$\sqrt{\ln|\ln(T-t)|}$ factor is
 the {\em nonlinear focusing Schr\"odinger equation} (the NLSE)
 with the critical exponent
  \be
  \label{1+}
   \tex{
\ii u_t=-\D u  - |u|^{p-1}u \inB \ren \times (0,T),
 \quad \mbox{where} \quad p=p_0=1+ \frac 4N,
 }
  \ee
  which is a fundamental model of water wave theory ($N=1$), nonlinear optics (mainly $N=2$),
   and plasma
  physics ($N=3$, as a limiting case of Zakharov's model of Langmuir waves, 1972).
 It seems that G.M.~Fraiman in 1985 \cite{Frai85}
  (see also \cite{Smi91}, and a full history of asymptotics
  in \cite[p.~115]{SulMon}) was the first,
 who formally derived the log-log correction to the blow-up rate
for the cubic case in dimension\footnote{There are two simple
misprints in \cite{Frai85} (in both Russian and English versions
that are identical) that are convenient to know: (i) the formula
(5.8) must read
 \be
 \label{corr1}
\tex{
 ``\ddot a=(1-|C_0|^2) P_{\rm cr}/a^3 \xi^2_{\rm eff}," \quad
   \big(\mbox{instead of} \quad  ``\ddot a=(1-|C_0|^2 P_{\rm cr})/a^3 \xi^2_{\rm eff},"
   \big).
   }
   \ee
In the process of its derivation, three lines above (5.8) the
reference ``(3.25)" (a nonexistent formula) seems should be
replaced by ``(3.17)".
 Then (5.14) implies (5.15), i.e.,
``$|C_0|^2-1 \sim 1/\ln^2 \t$", where $\t \sim -\ln(t_0-t)$, so
\ef{corr1} yields  the result:
 \be
  \label{corr2}
  \tex{
 \ddot a(t) \sim \frac 1{\ln^2 \t \, a^3(t)} \LongA a(t) \sim  \sqrt{\frac{t_0-t
 }{\ln|\ln(t_0-t)|}} \asA t \to t_0^-,
  }
  \ee
  with, according to (3.1), $\frac 1{a(t)}$  being the desired extra non-self-similar
   blow-up
  factor as in \ef{LL1} for $N=2$.

(ii) Thus, according to the correct asymptotics in \ef{corr2},  in
the final formula (5.16) in \cite{Frai85} (cf. p.~400 in the
Russian version), the second $\ln$ was incidentally missing, as
earlier (and for the first time?) noted in \cite[p.~116]{SulMon}.
Overall, Fraiman's derivation of log-log looks very solid and
formally well-justified from the point of view of standard
asymptotic methods; in particular, it does not use a
``homotopying" in the dimension $N \to 2^-$,
which was used in some other papers a few years later, until
perfect solution via solid operator-functional analysis approaches
developed by Perelman, Merle, Raphael, and others (see references
below)  in the twenty first century only.


 A single log-correction was suggested earlier by V.I.~Talanov in
the 1970s (see \cite{Vlas78}) and by D.~Wood in the 1980s  (see
also references in \cite{Land88}); V.E.~Zakharov also claimed
earlier derivation of the double log-log \cite{Zak91}, see some
details and results of the Russian School in \cite{Kos91}.
Self-focusing (blow-up) itself in the NLSE was proclaimed in the
middle of the 1960s, \cite{Chi64, Kell65, Ahm66}.}  $N=2$. An
alternative way of derivation of the log-log was presented in
\cite{Land88, Mes88} (1988).
 See the monograph \cite{SulMon} for a full history and extra
details, and \cite{Per01} for an extended list of further
qualitative results. These first results were based on formal
ideas saying that the $L^\iy$-blow-up rate for some classes of
solutions with $L^2$-norm slightly above that for the {\em ground
state} $Q(x)$, satisfies
 \be
 \label{LL1}
  \tex{
  \sup_x |u(x,t)| \sim \Big[\sqrt{ \frac {\ln|\ln(T-t)|}{T-t}}\,\, \Big]^{\frac N2}
  \asA t \to T^-.
  }
  \ee
 We refer to a number of more recent papers, where first
 qualitative and formal estimates obtained a proper mathematical
 justification \cite{Ang02, Budd01, FibM06, Fib07,
 Ken06, Li07, Mer04, Mer05, MerleR05, MerR052, Per01, Plan07, Rap06}.
    Many questions remain still open
 in view of the general complexity of the NLSE \ef{1+}, especially
 in the multi-dimensional geometry, with $N \ge 2$, where even the
 representation \ef{LL1} is questionable (other norms can be
 more appropriate).

 In general, blow-up such as \ef{LL1} for the NLSE \ef{1+}
 with the strong $L^2$-conservation property (the
 dynamical system being symplectic-Hamiltonian), can be
 characterized as corresponding to a {\em centre-subspace behaviour}
 relative to the 1D manifold of the rescaled ground states
 $\{Q_\o(x)\}$, though sometimes  this is again not well-understood and
   questionable (but this is of importance for other models with
 a different reason for the factor).

\subsection{THIRD $\sqrt{\bf log\,log}$: a
 cubic second-order reaction-diffusion model (1985)}

As rather customary in PDE theory, the above  log-log factor was
almost simultaneously (with Fraiman's research \cite{Frai85})
detected in reaction-diffusion theory, and was
 involved  into mathematical blow-up area by Friedman and McLeod in
1986 \cite{FM86}, who studied the following parabolic equation
from plasma physics:
 \be
 \label{par1}
  \tex{
   u_t=  u^2(u_{xx}+u) \inB (-L,L) \times (0,T),
 \quad u(\pm L,t)=0 \whereA
    L > \frac
   \pi 2.
   }
   \ee
    More precisely, according to B.C.~Low
     \cite{Low73}\footnote{``The plasma and current sheet collapse
     is in fact forming a singularity, suggested by B.C.~Low...",
     \cite[p.~3]{Vain07}.},
  \ef{par1} can be considered as a model of resistive diffusion of a
 force-free  magnetic field in a plasma.
 The multi-dimensional version of the model assumes replacing
 $u_{xx}$ by the Laplacian $\D u$. Actually, Low proposed the original
 nonstationary model, in which a contraction of the magnetic field
 becomes infinite after a finite time provided that the
 free-force condition is always satisfied.

    It is worth mentioning here that plasma physics with magnetic fields is a permanent
 source of various reaction-diffusion problems, where evolution of ``frozen"
 magnetic field is typically described by fast-diffusion-like and
 other quadratic ambipolar
 mechanisms. For convenience and better understanding a hierarchy of the corresponding parabolic models,
   we present a typical multi-term RD-type equation for the
 magnetic field $\BB$ \cite{Bra94}:
  $
  \BB_t= \n \times ({\bf v}_{\rm n} \times \BB + {\bf v}_{\rm D}
  \times \BB - \l \n \times  \BB),
  $
  where $\l$ is the ordinary
  magnetic diffusivity, and
 ${\bf v}_{\rm D}={\bf v}_{\rm i}- {\bf v}_{\rm n}$ is the
  drift velocity between ions and neutrals, which is proportional
to the Lorentz force,
 $
{\bf v}_{\rm D} \sim \frac {(\n \times \BB)\times \BB }{4 \pi
\rho_{\rm i} \nu_{\rm in}}
 $
 (this is the mechanism of cubic nonlinearities).
 Overall, this yields the following equation:
 $$
 \tex{
  \BB_t= \n \times \big[{\bf v}_{\rm n} \times \BB +
\frac {(\n \times \BB)\cdot \BB} {4 \pi \rho_{\rm i} \nu_{\rm
in}}\, \BB - (\l+\l_{\rm AD}) \n \times \BB \big] \whereA \l_{\rm
 AD}=
 \frac{\BB^2}{4 \pi \rho_{\rm i}\nu_{\rm in}}
 }
 $$
 is the ambipolar diffusion coefficient.
  In  addition,
 there can appear
  the related 1D  complex cubic PDEs of the PME-type (i.e., unlike \ef{par1}, with a
  divergent diffusion operator)
   $
   {\mathcal B}_t=\big(\frac {\mathcal
 B}2(|{\mathcal B}|^2)_z\big)_z
  $
   (without ordinary magnetic diffusion, no neutral velocity, and very small
  ionization fractions), which creates sharp structures in finite
  time, \cite[p.~L92]{Bra94}.

  First qualitative results for
 \ef{par1} including conjecture of blow-up were obtained in
   \cite{Wat85,
 Wat86}\footnote{The collapse ``..has been predicted analytically by Watterson (1986),
  who used a simplified one-dimensional, force-free model to show that a radially decreasing
   magnetic
 field profile with a reversal of the axial field cannot be maintained indefinitely",
  \cite[p.~508]{Ste87}.}.
 On the other hand, \ef{par1} appears in curve-shortening flows by
 curvature\cite{Ang91}. Namely, the equation $X_t= k N$ of the
 evolution of a curve $X:S^1\times[0,T) \to \re^2$ on the plane,
 where $k$ is the curvature and $N$ is the unit normal, is reduced
 to
  $$
  k_t=k^2(k_{\th\th}+k) \quad \mbox{on} \quad S^1 \times (0,T)
  \quad\big(S^1=\{0 \le \th<2\pi\} \,\,\,\mbox{is the unit circle in $\re^2$}\big).
   $$
 On derivation,  history, blow-up patterns, and other details, see  papers \cite{Ang91, AngVel1}
  to be referred to again. For earlier results and blow-up
  similarity solutions for such a flow, see e.g., \cite{Arb86}.

First, Friedman and McLeod \cite[\S~2]{FM86} proved blow-up of
solutions if $L
> \frac \pi 2$ by using standard at that time barrier and Maximum Principle
techniques. Among more delicate results, they proved a rather rare
(again at that time) result on the {\em blow-up set} of the
solution $u(x,t)$:
 \be
 \label{Bl1}
 B[u_0]=\big\{x \in (-L,L): \quad \exists \,\,\{x_n\} \to x,
 \,\,\{t_n\} \to T^- \,\, \mbox{such that}\,\, u(x_n,t_n) \to
 +\iy\big\}.
 \ee
By the definition, $B[u_0]$ is closed. Then,  for smooth symmetric
monotone data \cite[\S~4]{FM86}
 \be
 \label{Bl2}
  \tex{
  u_0(-x) \equiv u_0(x), \quad u_0'(x) \le 0 \,\,\,
  \mbox{for}\,\,\, x \ge 0
  \LongA B[u_0]=[-\frac \pi 2, \frac \pi 2].
   }
   \ee

   The next question is then the description of the blow-up
   behaviour in $B[u_0]$ as $t \to T^-$, which turns out to be a
   difficult issue.
 Namely, the authors showed by some involved
mathematics associated with the equivalent integral equation,
obtained by dividing by $u^2$ and inverting the operator $D_x^2
+I$ on the right-hand side, that the blow-up behaviour does not
obey the dimensional similarity law associated with the separation
of variables:
 \be
 \label{par2}
  \tex{
  u(x,t) \sim \frac 1{\sqrt{T-t}}\, \th(x) \asA t \to T^-
  \LongA \frac 12\, \th=\th^2(\th''+\th).
   }
   \ee
More precisely, unlike \ef{par2}, they showed that the behaviour
 is governed by the non self-similar law
 \be
 \label{par3}
 \tex{
  u(x,t) \sim \frac {A(t)}{\sqrt{T-t}}\, \cos x
  \whereA A(t) \to +\iy \asA t \to T^-.
   }
   \ee
In addition, it was mentioned on the last page \cite[p.~80]{FM86}
 that in Watterson's Thesis in 1985\footnote{I.e., in the same 1985 year as Fraiman's
 conclusion for the cubic NLSE \cite{Frai85}! What nonlinear theory can explain this historical
  ``time-evolution" coincidence?}
   \cite{Wat85} the \ef{AA1} ansatz was
 made\footnote{As explained, ``... on the basis of numerical evidence",  which is
 hardly believed even in the twenty-first century: log-log in
 blow-up is extremely difficult to catch numerically since a huge solution growth of hundreds of orders
  is necessary.}.
An involved proof of \ef{AA1} and the behaviour \ef{par3} was
achieved later on in \cite{AngVel1} with a high sharpness:
 the result therein was
 \be
    \label{par4S}
 \tex{
 A(t) = \sqrt{\ln|\ln(T-t)|}(1+o(1)) \asA t \to T^-
 }
 \ee
 (we will explain why the multiplier $1 \cdot \sqrt{...}$ always occurs all the time here on matching).

It is convenient for further use to discuss the origin of the
structure \ef{par3}.
 The $\cos x$ is obviously a formal stationary solution of \ef{par1},
  \be
  \label{par31}
  f(x)=\cos x: \quad f_{xx}+f=0,
 \ee
 which nevertheless does not satisfy the Dirichlet boundary
 conditions in \ef{par1}, since $L> \frac \pi 2$. Therefore,
 \ef{par3} can be interpreted as an evolution as $t \to T^-$ along the 1D
 manifold $\{\mu \cos x, \,\, \mu>0\}$ of stationary solutions,
 where the resulting slow-growing factor \ef{AA1} represents a
 kind of centre manifold (or subspace) behaviour corresponding to
 the linearized rescaled operator. Of course,  the crucial
 argument establishing \ef{AA1} is fully dependent on the
 matching of the blow-up structure \ef{par3} with bounded
 solutions near $x=\pm L$ through the special behaviour at the
 localization end points $x = \pm \frac \pi 2$.

\subsection{On a ``topological" matching: main idea}

For further applications to more difficult nonlinear PDEs, we need
 to explain the qualitative origin of our rather ``topological"
 asymptotic analysis. We claim that the actual origin of
the log-log factor \ef{AA1} in the RD model \ef{par1}
 can be associated
with the existence of the log-travelling wave solutions (log-TW)
of the equation \ef{par1} of the following
 form
 \be
 \label{par6}
  \tex{
 u_{\rm log}(x,t)=\frac 1{\sqrt{T-t}}\,  g(\eta), \quad \eta= x+ \l \ln(T-t)
 \quad (\l \not = 0).
  }
  \ee
Substituting \ef{par6} into \ef{par1} yields the following ODE:
 \be
 \label{par7}
  \tex{
 \frac 12 \, g - \l g'= g^2(g''+g).
  }
   \ee
 These solutions can be obtained as the result of the invariance of
 \ef{par1} relative to a one-parameter group of scaling-translating
 transformations; this was first proved by Ovsiannikov in 1959
 \cite{Ov59}.

We show that  similarity solutions \ef{par6} play a key role at
the end points $x= \pm \frac \pi 2$ of the localization domain $x
\in ( -\frac \pi 2, \frac \pi 2)$, and serve  as a ``transitional
mechanism" from the Inner Region-II of bounded solutions for $|x|
> \frac \pi 2$ into the internal Inner Region-I $\{|x|< \frac \pi
2\}$ with the blow-up behaviour \ef{par3}.
 This idea is explained in Figure \ref{F1}, which will be used in
 greater detail below.
Since in some asymptotic $\eta$-interval, solutions $g(\eta)$ of
\ef{par7} have an extra logarithmic correction in the asymptotics,
i.e.,
 \be
 \label{par8}
 g(\eta) \sim (-\eta) \sqrt{\ln(-\eta)} \quad \mbox{for}
 \quad \eta \ll -1
 \ee
(in this intermediate region, the term $g^3$ in \ef{par7} ought to
be neglected), the combination of two logs:


 (i) $\sqrt{\ln(-\eta)}$ in \ef{par8}, and


   (ii)  $\ln(T-t)$ in the
$\eta$-variable in \ef{par6},


 \noi lead, on a special formal matching, to the double logs
in \ef{AA1}  in the RD model \ef{par1}\footnote{This idea was
developed by the author  in  discussions with Herrero  and
Vel\'azquez in
 Dpto. de Matem\'atica Aplicada, Unversidad Complutense de
Madrid, 1992, as explains in a few lines in \cite[p.~308]{SGKM}.}.
The later difficult proof such a blow behaviour with the factor
\ef{AA1} in Angenent--Vel\'azquez \cite{AngVel1} of the behaviour
\ef{par3}, \ef{AA1} essentially uses the rescaled variables
corresponding to the log-TWs moving frame \ef{par6} and other
related issues.

In what follows, to avoid very complicated technical calculus and
further justification (which are actually nonexistent for
higher-order PDEs under consideration), we prefer to keep the
ideology of the above ``topological" matching with a clear
geometrical meaning to be explained in Figure \ref{F1} in Section
\ref{S2}.

\subsection{No log-log for the divergent RD model with finite propagation}

The corresponding cubic RD model with the divergent diffusion
operator has the form
 \be
 \label{do1}
 u_t=(u^3)_{xx} + u^3 \inB \re \times (0,T).
  \ee
Answering the question \ef{Q1}, the generic blow-up of nonnegative
solutions is then
 described by  the following  explicit {\em Zmitrenko--Kurdyumov
solution} \cite{SZKM2}:
 \be
 \label{RD.31}
  \tex{
 u_{\rm S}(x,t)= \frac 1{\sqrt{T-t}}\, \th(x)
  \whereA
 \frac 12  \, \th= (\th^3)'' +  \th^{3} \inB \re, \quad \mbox{so}
  }
  \ee
 \be
 \label{RD.4}
 \mbox{$
   \th(x) = \left\{ \begin{matrix}  \frac {\sqrt 3}2  \,\cos \big(\frac x3\big)
   \quad \mbox{for} \quad |x| \le  \frac {3\pi}2, \smallskip\\ \,\,
   \,\,\,\,\quad 0\quad \qquad \mbox{for} \quad |x| >
   \frac {3\pi}2.
 \end{matrix}
    \right. $}
   \ee
 Concerning the asymptotic stability and other generic evolution properties  of the particular solution \ef{RD.31},
 see \cite[Ch.~4]{SGKM} and references therein.
In other words, the weak solution \ef{RD.31}, \ef{RD.4} shows that
the blow-up in localized in the domain $(- \frac{3\pi}2,
\frac{3\pi}2)$, which the heat does not penetrate from for all $t
\in (0,T)$. The rate of blow-up in \ef{RD.31} is purely
self-similar and does not contain any extra log-log or other
factors as in \ef{AA1}.

Indeed, returning again to \ef{Q1}, this absence of the log-log
factor can be directly connected with the solvability of the ODE
in \ef{RD.31}, which admits a good weak compactly supported
profile $\th(x)$ given in \ef{RD.4}. This profile is suitable for
both the Cauchy problem in $\re \times (0,T)$ and for the IBVP in
$(-L,L)\times (0,T)$ provided that
 \be
 \label{p2}
 \tex{
 L \ge \frac {3\pi}2.
  }
  \ee
 If \ef{p2} is violated, other boundary conditions at the end
 points $x= \pm \frac {3\pi}2$ can be imagined that can lead to
 extra blow-up factors. Though such posed IBVP for \ef{do1} are
 typically looking rather artificial unlike the above and given below more natural
 parabolic models.

\subsection{Main models and results: explaining log-log in higher-order RD
models via transition by log-TWs}

 We claim that
the same blow-up factor as in \ef{AA1} exhibits more universality
and occurs in other nonlinear scaling-invariant cubic parabolic
models of higher-order. Moreover, we also claim that, similar to
\ef{par6},
 $$
  \fbox{$
\mbox{$\sqrt{\log\,\log}$ is generated by the log-TWs transition
mechanism.}
 $}
 $$

 As a key example, we use the
fourth-order RD equation, as a natural extension of \ef{par1},
with zero Dirichlet boundary conditions:
 \be
 \label{p1}
  \left\{
  \begin{matrix}
 u_t= u^2(-u_{xxxx} + u) \inB (-L,L) \times (0,T), \ssk\\
  u=u_x=0
 \,\,\, \mbox{at} \,\,\, x= \pm L \quad (L > L_0).
 \qquad\quad \,\,\,
  \end{matrix}
  \right.
  \ee
 Of course, the problem for this higher-order parabolic equation
loses any traces of order-preserving, comparison, the Maximum
Principle, and barrier features that  are always key ingredients
in the mathematical study of the second-order parabolic equations.
Moreover, even existence-uniqueness results (local) for the
degenerate equations such as \ef{p1} are not properly settled.
Nevertheless, we put a blind eye on those difficulties, and
concentrate on predicting the log-log factor in the regional
blow-up behaviour. Actually, there is some space to avoid such
local difficulties: for instance, we can consider uniformly
strictly positive solutions with such data on the boundary that
the similarity blow-up such as \ef{RD.31} is impossible (but then
some extra speculations are necessary). For positive solutions
$u(x,t) \ge \d_0>0$, classic parabolic theory \cite{EidSys, Fr}
applies to guarantee local existence and uniqueness of smooth
(moreover, analytic) solutions; see a comment at the beginning of
Section \ref{S2} concerning other weaker solutions $u \ge 0$.

 Any rigorous proof of such a blow-up log-log behaviour for \ef{p1} then becomes rather
 illusive (e.g., more difficult than for the NLSE \ef{1+}, and any second-order RD-type problem
  considered before), while
 formal using of the log-TWs \ef{par6} is the only source of a
 proper justification.

 In Section \ref{S2} we explain the blow-up log-log phenomenon for
 the model \ef{p1}.
In Section \ref{S2m3}, we briefly explain how the same factor
occurs in the sixth-order parabolic problem
 \be
 \label{m3}
 u_t=u^2(u_{xxxxxx}+u),
  \quad u=u_x=u_{xx}=0
 \,\,\, \mbox{at} \,\,\, x= \pm L \quad (L > L_0).
  \ee
 In Sections \ref{S3} and \ref{S4}, we show that blow-up is simply
 self-similar  for the better divergent models: the {\em fourth-order
 porous-medium equation with source} (the PME--4),
  \be
  \label{p3}
  u_t= -(u^3)_{xxxx}+u^3,
  \ee
 and for the {\em thin film equation with source} (the TFE--4)
  \be
  \label{p4}
  u_t= -(u^2 u_{xxx})_{x}+u^3.
  \ee
  Without loss of generality,
for both PDEs \ef{p3} and \ef{p4},
 we consider the Cauchy problem
in $\re \times (0,T)$ or the IBVP as in \ef{p1}.

\subsection{On extensions to wave and nonlinear dispersion
equations}

Finally, in Section \ref{S5}, we show that a similar blow-up
log-TW mechanism  makes it possible to reconstruct log-log factors
for some blow-up patterns for other nonlinear PDEs: for the {\em
quasilinear wave equation} (the QWE--4)
 \be
 \label{QWE1}
 u_{tt}=u^2(-u_{xxxx}+u),
  \ee
  and for the {\em nonlinear dispersion equation} (the NDE--3)
 \be
 \label{NDE1}
 u_{t}=u^{3}(u_{xxx}+u) \quad (\mbox{note ``$u^{3}$" instead of
 ``$u^2$"}).
  \ee
 Actually, the formal matching does not change at all in
 comparison with the parabolic case. Of course, equations
 \ef{QWE1} and \ef{NDE1} contain other singularity phenomena such
 as, in view of nonlinear dispersion mechanisms involved (the local speed of propagation
 depends on the value of $u$ itself), formation of various shocks, with rather difficult and not
 fully justified ``entropy-like"  mathematics (see \cite{GalNDE5} as a guide for
 higher-order NDEs and \cite{GalQWE} for QWEs--4). We do not take these phenomena into account,
 and just formally show that there exist some particular blow-up patterns with
 log-log property, and do not discuss any of their structural
 stability properties, which can lead to rather obscure
 mathematics.





\subsection{Towards consistency of ``topological" matching}
 \label{S1N}

 As a final comment, the author  emphasizes that the
 non-rigorous and rather rough geometric  nature of the presented results makes no problem for
 him, since, according to his almost thirty years experience in
 proving various blow-up results, for some of the above
 higher-order models\footnote{
 ``The main goal of a mathematician is not
 proving a theorem, but an effective investigation of the
 problem..." (A.N.~Kolmogorov, 1980s; the author apologizes for a
non-literal translation from the Russian).},
  \be
  \label{DD1}
   \fbox{$
   \begin{matrix}
 \mbox{proving the proposed blow-up $\sqrt{\ln\,\ln}$
 asymptotics} \ssk\\
 \mbox{is not possible in a reasonable finite time.}
  \end{matrix}
   $}
    \ee
    In fact,
 it is not that easy even to explain how difficult those problems
 are. In  Section \ref{S2.5} we present some operator theory
 comments, and meantime
clarify the above point of view, partially  expressed in \ef{DD1},
in a more formal but consistent manner, as follows:


 (i) It should be noted that the higher-order degenerate
singular PDEs such as \ef{p1}, \ef{m3}, \ef{p4}, etc., though
having already good and clear applications, have quite {\em
obscure} local existence-uniqueness-entropy-... theory. The
problem is not around blow-up (this can be done rigorously for
some classes of solutions), but even small solutions are extremely
oscillatory; see similar examples in \cite{Gl4} and
\cite[Ch.~3-5]{GSVR}. To say nothing about the hyperbolic problem
\ef{QWE1}, which cannot have a smooth local solution at all, since
in view of an obvious nonlinear dispersion mechanism, shock (and
rarefaction) waves can appear at any suitable point.
 Note that by no means \ef{QWE1} is a (strictly) hyperbolic system
 with good nowadays entropy theory. Moreover,
 it can be expected that PDEs such as \ef{QWE1}
 principally cannot have a local uniqueness and well-posedness
 theory.
 Similar
phenomena of singularity and non-uniqueness exist for \ef{NDE1}.


 (ii) Overall, in general, in view of a complete absence  of
local existence-uniqueness and entropy theory, studying these
higher-order models, we actually deal with {\em huge bundles}
(so-called, {\em flows} in dynamical system theory) of solutions,
which can be responsible for various problems settings, from the
Cauchy problem up to infinitely many free-boundary problems
(FBPs). This is a typical, and unavoidable, feature of modern PDE
theory: for equations with non-monotone, non-coercive,
non-variational... operators of higher-order, proving existence,
uniqueness, entropy, ... results along the lines great classic
theory from the twenties century is not only very difficult, but
can be impossible in principle.


 (iii) Therefore, dealing with bundles (flows) of blow-up
solutions, we propose a rather rough asymptotic method, which, as
must be admitted,  does not specify many particular features of
solutions involved, but is able to detect some universal property
of $\log \log$-factor, which holds {\em for all such blow-up ones,
regardless which problems (Cauchy, FBPs, Neumann, Robin, Florin,
Stefan,...) are posed.}
 Further refining  of the asymptotic method, will require a
sharper posing of the Cauchy or FBPs, i.e., matching with boundary
behaviour, which is completely different (e.g., extremely
oscillatory) for various problems.
 Inevitably, one faces  difficult local
 existence/uniqueness/entropy/etc. aspects that are unclear for
 the most of models, which cannot be a subject of a  single
 paper.


(iv) In a natural sense, the proposed ({\em ``topological"}, but
very rough, if you wish) matching well corresponds to common
understanding coming from the theory of asymptotic series, which
are not converging but first terms correctly describe the
behaviour of the functions relative to a small parameter $\e \to
0$. From planetary motion study in the seventeenth century
(Newton-Halley's period), it is known that taking many terms of
asymptotic expansion gives worse results, while first terms  are
sufficient to predict the behaviour. In a certain  sense,  a
similar asymptotic phenomenon happens in our study: trying to
improve the expansion, we involve extra, highly oscillatory and
singular terms, which  are deeper connected with the whole
solution bundle that  we are not  unaware of.


(v) Thus, even in the clear presence of a certain fear of a
justified criticism from the attentive Readers, the author prefers
to keep the style of the asymptotic analysis more based on a {\em
topological matching}, rather than a traditional metric one, where
all the terms are carefully estimated to convince. Of course, a
standard  balancing of all the differential terms  can be done for
all the higher-order models, but not more than that: any further
refining the expansion enters the oscillatory-singular solution
bundles and any meaning and significance   of the blow-up
expansion will be  lost and become illusive.

 \section{How log-TWs imply the log-log blow-up factor}
 \label{S2}

Thus, we fix the model \ef{p1} as the basic one for explaining the
occurrence of the log-log factor \ef{AA1} in the generic blow-up.
Concerning questions of existence and uniqueness of solutions of
such higher-order equations, we refer to the pioneering paper by
Bernis and Friedman \cite{BF1}, where a general approach to
constructing such nonnegative solutions has been developed.
 Clearly, their methods of ``singular parabolic $\e$-regularizations"
 for constructing non-negative solutions of the Cauchy problem or
 free-boundary problems (depending on parameters)
apply to a very wide class of equations including \ef{p1},
\ef{m3}, and many others, linear or nonlinear.
 On the
other hand, for a number of higher-order parabolic equations,
proper solutions of the Cauchy problem can be oscillatory and of
changing sign (see \cite{BMcL91, Gl4},  \cite[p.~152]{GSVR} and
references therein), but these local interface aspects can be
neglected
 at this stage
 while
dealing with peculiar properties of blow-up of very large
solutions.

\subsection{On nonexistence of blow-up separable-variable solutions}

Bearing in mind \ef{Q1}, this is the first crucial step of any
blow-up study. Thus, we look for the standard self-similar blow-up
for \ef{p1} in the form as in \ef{RD.31},
 \be
 \label{p41}
  \tex{
 u_{\rm S}(x,t)= \frac 1{\sqrt{T-t}}\, \th(x)
  \whereA
 \frac 12  \, \th= \th^2(-\th^{(4)} +  \th).
  }
  \ee
  Recall that we consider either the Cauchy problem for the equation in \ef{p1},
  so we demand a localized patter satisfying
   \be
   \label{p5}
    \th(x) \to 0 \asA x \to \pm \iy,
     \ee
or the IBVP as in \ef{p1}, i.e., one needs a proper solution
$\th(x)$ such that
 \be
 \label{p6}
 \th=\th'=0 \quad \mbox{at} \quad x= \pm L.
  \ee

It is easy to see that such solutions of the ODE in \ef{p41} do
not exist. Dividing the ODE by $\th^2$, multiplying by $\th'$, and
integrating by parts yields
 \be
 \label{p7}
 \tex{
 \frac 12\, \ln|\th|=-\th''' \th +\frac 12 \, (\th'')^2+ \frac 12\, \th^2+C,
 }
 \ee
 where $C \in\re$ is a constant of integration. This identity is
 true in the standard classic sense for smooth positive solutions
 $\th(x)>0$, and also remains valid for a proper class of
 solutions of changing sign having sufficiently regular zeros
 (this is a question of a suitable functional setting, not
 accented here). We fix this negative result:

 \begin{proposition}
 \label{Pr.Non}

 The identity $\ef{p7}$ does not allow nontrivial sufficiently smooth, bounded, and
  vanishing solutions $\th(x)$ satisfying  $\ef{p5}$ or $\ef{p6}$.

  \end{proposition}

This negative result can be associated with some king of a strong
local monotonicity of a nonlinear elliptic operator induced by the
equation in \ef{p41}.

 \subsection{On self-similar blow-up for the Neumann or periodic problem}

 Note that, for the Neumann conditions for the IBVP,
 \be
 \label{p8}
 \th''=\th'''=0 \quad \mbox{at} \quad x= \pm L,
  \ee
the ODE in \ef{p41} can admit a good proper solution $\th(x)$, and
then the blow-up becomes purely self-similar as in \ef{p41}.
 The same can occur for the periodic problem:
 \be
 \label{p8p}
 \th(-L)=\th(L), \quad \th'(-L)=\th'(L), \quad\th''(-L)=\th''(L),
 \quad\th'''(-L)=\th'''(L).
  \ee
 We arrive at the problem of  asymptotic structural stability of $th$ as $t \to T^-$,
which can be solved in a standard manner by considering the {\em
rescaled equation} for the rescaled solution:
 \be
 \label{p9}
  \tex{
  v(x,\t)= \sqrt{T-t} \, u(x,t), \,\,\, \t=- \ln(T-t) \LongA
 \fbox{$
v_\t = v^2(-v_{xxxx}+v)- \frac 12\, v.
  $}
 }
   \ee
It is key that \ef{p9} is a gradient dynamical system with the
Lyapunov function
 \be
 \label{10}
  \tex{
  \frac {\mathrm d}{{\mathrm d}\t}\big[\frac 12 \int(v_{xx})^2
   -\frac 12 \int v^2 + \frac 12 \int \ln |v| \big]=- \int
   \frac{(v_\t)^2}{v^2} \le 0.
    }
    \ee
 As in standard blow-up theory, the main difficulty in proving the
 stabilization as $\t \to +\iy$ in the rescaled PDE \ef{p9} will
 be {\em a priori} bounds on the orbit, i.e., that
  \be
  \label{11}
   \tex{
  |v(x,\t)| \le C \andA \sup_x |v(x,\t)| \not \to 0 \asA \t \to +
  \iy.
  }
   \ee
 Then, if the set of solution of the ODE \ef{p41}, \ef{p8} or \ef{p8p} is
 discrete, then there exists a unique stationary profile $\th(x)$
 such that
  \be
   \label{12}
    v(x,\t) \to \th(x) \asA \t \to +\iy,
    \ee
 in a suitable metric (e.g., uniformly for strictly positive and
 hence classic solutions).

Finally, let us  notice that the solvability and multiplicity for
the Neumann problem \ef{p8} or a periodic one \ef{p8p}  are not
easy (though the variational structure suggests existence of many
periodic solutions), essentially depends on the length parameter
$L$, and can deliver some striking  surprises.

\subsection{Stationary profiles and a bound on $L$}

Thus, assuming that a standard separable blow-up solution does not
exist, we then ought to find the stationary one (cf. \ef{par31}),
 \be
 \label{13}
 f(x): \quad -f^{(4)}+f=0 \quad \mbox{on} \quad (-L_0,L_0), \quad
 f=f'=0 \quad \mbox{at} \quad x= \pm L_0.
  \ee
Looking for an even profile, we set ($L$ denotes $L_0$)
 \be
 \label{14}
 f(x)= \cos x + C \cosh x \LongA \left\{
  \begin{matrix}
  \cos L=-C \cosh L, \ssk\\
  \sin L=C\sinh L,\,\,\,\,
   \end{matrix}
   \right.
   \quad \mbox{or} \quad \tan L = - \tanh L.
   \ee
 The last algebraic equation yields the first positive root $L_0
 \in (\frac \pi 2, \pi)$, which is taken into account in the statement
 \ef{p1} in order to avoid existence of classic stationary
 solutions and hence no blow-up.

 \subsection{Formal matching with log-TW behaviour: a ``topology" approach}

\begin{figure}
\centering
\includegraphics[scale=0.85]{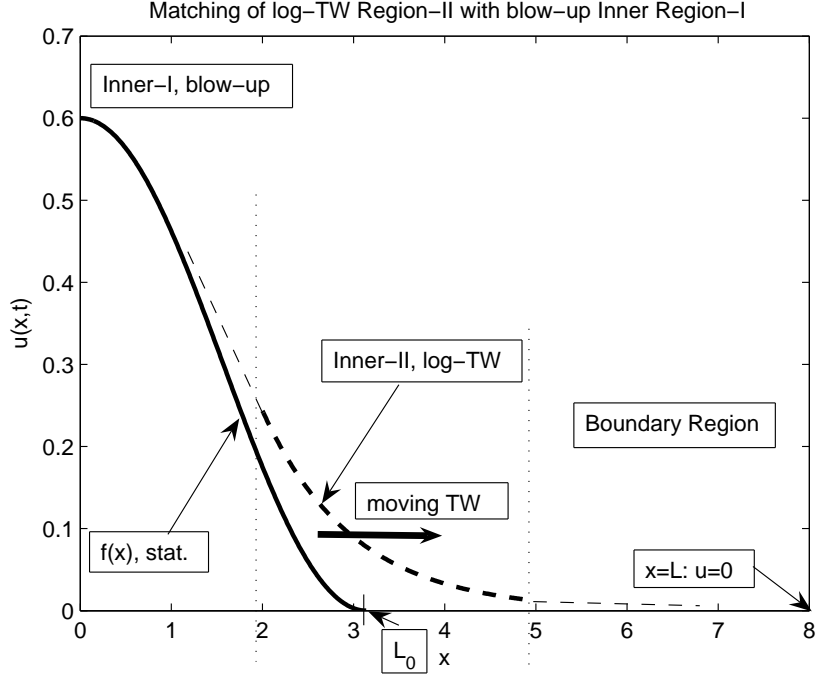} 
\vskip -.3cm \caption{\small On the idea of a schematic
``topological"  matching of two Inner Regions.}
   \vskip -.3cm
 \label{F1}
\end{figure}

 Thus, assuming that $L>L_0$, similar to \ef{par3}, we suggest that  the
 blow-up evolution is close to the manifold of stationary solutions,
 i.e., we study the solutions with the behaviour
 \be
 \label{15}
  \tex{
   u(x,t) \sim \frac{A(t)}{\sqrt{T-t}}\, f(x),
    }
    \ee
 where $f(x)$ is a stationary solution given by \ef{14} normalized by its value
 at the origin,
  \be
  \label{16}
   f(0)=1+C \quad (C>0).
    \ee
The further strategy of matching is explained in Figure \ref{F1},
which suggests to match the Inner Region-I with the strong and
fast blow-up asymptotics \ef{15} with the log-TW Inner Region-II
situated close to the end point $x=L_0$ (we consider the
$x$-symmetry geometry).

We are assuming that in Region-II, the solutions are sufficiently
small, so we take into account  the diffusion-like operator only,
 \be
 \label{17}
  u_t= - u^2 u_{xxxx} \quad (x \approx L_0).
  \ee
 Then, taking the log-TW \ef{par6} yields the following ODE:
 \be
 \label{18}
  \tex{
  \frac 12\, g - \l g'= -g^2 g^{(4)}.
 }
   \ee
 It is now key that the ODE \ef{18} admits the following
 asymptotics:
  \be
  \label{19}
   \tex{
   g(y)= (-y)^2 \sqrt{\ln(-y)}(1+o(1)) \asA y \to - \iy.
   }
   \ee
Actually, the leading quadratic behaviour $\sim y^2$ well
corresponds to the quadratic ``boundary layer" with the stationary
structure
 \be
 \label{20}
 f(x) = C_1(L_0-x)^2(1+o(1)) \asA x \to L_0^- \quad(C_1>0).
  \ee
The logarithmic correction in \ef{19} has the same type as in
\ef{par8}, so it is supposed to generate the same log-log factor
as in \ef{AA1}.

More precisely, according to \ef{par6}, we need to match the
behaviour in Region-II,
 \be
 \label{21}
  \tex{
 u(x,t) \sim \frac 1{\sqrt{T-t}}\, \big[-x-\l \ln(T-t)\big]^2\,
 \sqrt{\ln(-x-\l \ln(T-t))},
 }
  \ee
 with that in blow-up Region-I, where, by the assumption \ef{15},
 the following holds:
  \be
  \label{22}
   \tex{
  u(x,t) \sim \frac {A(t)}{\sqrt{T-t}}\, C_1(x-L_0)^2 \quad \mbox{for} \quad
   x \approx L_0^-.
  }
   \ee
 Note that the expansion \ef{21}
 should be treated as already the {\em extended} one from
the corresponding asymptotic region. A more detailed expansion for
smaller values of $u$ will be inevitably affected by the singular
and oscillatory behaviour of solutions close to $x=L_0$ and $x=L$,
which we would like to avoid at this stage and refer to the
discussion in Section \ref{S1N}.  More definitely, the expansion
\ef{21} up to some admissible perturbations (and partially
\ef{22}, which is more regular) precisely defines the required
class of solutions of the IBVP or any of FBP problems. Of course,
there are other solutions of the PDE under consideration with a
completely different blow-up behaviour having nothing to do with a
double-log factors.


We omit at this moment the necessary (delicate and questionable)
stage in matching, where a suitable time-dependence of the
log-velocity $\l=\l(t)$ is assumed, \cite[p.~308]{SGKM}.
Currently, as usual in matching concepts and theory, two manifolds
of solutions such as \ef{21} and \ef{22} admit matching if the
corresponding pairs of leading multipliers are overlapping at some
intermediate values of $x$, possibly, different for both pares.
 One can object that the log-TW expansion in \ef{21} describes
 {\em moving} wave (see the arrow in Figure \ref{F1}),
  unlike the {\em standing} one in \ef{22},
 which, at the first sight, makes matching rather suspicious.
 We should mention in connection with that
 the log-TW in Region-II moves with the logarithmic
 speed $\sim |\ln(T-t)|$, i.e., very slow in comparison with the power rate
 $\sim \frac 1{\sqrt{T-t}}$
 of the blow-up divergence in Region-I. In
other words, relatively,  this can be classified as an {\em
effectively standing} wave, which makes the matching correct and
possible after necessary extra arguments.

 Again carefully looking at the manifolds \ef{21} and \ef{22}, we
 see that the two pairs of different multipliers that are, as $t
 \to T^-$,
 \be
 \label{23}
   \fbox{$
(-x-\l \ln(T-t))^2\,\,\&\,\,(x-L_0)^2 $}
 \,\,\, \mbox{and} \,\,\,
  \fbox{$
\frac 1{\sqrt{T-t}}\,
 \sqrt{\ln(-x-\l \ln(T-t))}\,\,\&\,\, \frac {A(t)}{\sqrt{T-t}},
 $}
  \ee
 admit natural ``structural matching"
 on some compact subsets
 (different for both) in $x$.
Actually, according to typical concepts of asymptotic analysis,
one needs to observe and to get ``similar geometric forms" in
space and time of matched flows. This can be definitely done for
the pairs in \ef{23}, where the first pair of space distributions
are both similar and quadratic, while the second pair assumes to
be similar for some special functions $A(t)$ only. We must admit
again that by no means our goal is to prove this matching
rigorously (whatever ``proof" means).

  Then setting $x \approx L_0$ in the
 second pair leads to the conclusion \ef{AA1}: as $t \to T^-$,
  \be
  \label{24}
  \tex{
  \fbox{$
  \frac 1{\sqrt{T-t}}\,
\sqrt{\ln(-\l \ln(T-t))} \sim \frac {A(t)}{\sqrt{T-t}}
   $}
 \LongA A(t) \sim \sqrt{\ln|\ln(T-t)|}(1+o(1)).
  }
  \ee
Note that $\l$ even being a slow growing function completely
disappears in this matching, so even this rough accuracy
guarantees the coefficient ``1" (and $\l=\l(t)$ is contained in
the $o(1)$-term only, which is not of importance).

\subsection{On some standard metric estimates and operator balance:
the origin of principal difficulties}
 \label{S2.5}

Let us return to a more mathematical (rather than geometrical))
PDE meaning of these speculations. In fact, this suffices to
balance differential terms and perturbations in the rescaled
equation for $v$ in \ef{p9},
  where $\t=-\ln(T-t) \to +\iy$ as $t \to T^-$.
Actually, Figure \ref{F1} is more suitable for the rescaled
function $v(x,\t)$, which is expected to have a slow logarithmic
growth as $\t \to + \iy$ unlike the fast  algebraic one for
$u(x,t)$ as $t \to T^-$. We next need to rescale as in \ef{p9} the
asymptotic flows \ef{21} and \ef{22}, etc.

We begin with a standard approach to such a matching by starting
with the inner R-I rescaled expansion \ef{21}, i.e.,  $v(x,\t)
\sim A(\t) f(x)$,  which naturally require introducing the new
variable $w$ and the corresponding rescaled equation:
 \be
 \label{22N}
  \left\{
 \begin{matrix}
 v(x,\t)=A(\t) w(x,\t) \andA A^2(\t) \,{\mathrm d}\t= {\mathrm
 d}s, \ssk\\
  w_s= w^2(-w_{xxxx}+w)+ \frac 1{2A^2(\t)}\, w -
 \frac{A'(\t)}{A^3(\t)}\, w.\,\,
  \end{matrix}
   \right.
  \ee
Thus, here $A(\t)$ is unknown, but \ef{AA1} suggests that
 \be
 \label{A1}
  \tex{
 A(\t) \sim \sqrt{\ln \t} \andA s \sim \t \ln \t \, \,\, \big(\t \sim \frac{s}{\ln s}\big) \asA \t \to + \iy
  \quad (\mbox{expectation}).
  }
  \ee
Indeed, the perturbed equation \ef{22N} with the leading (by
\ef{AA1}) perturbation $\sim O\big(\frac 1{A^2(\t)})$ describes
the stabilization to the stationary profile:
 \be
 \label{23N}
 w(x,\t) \to f(x) \asA \t \to +\iy \whereA f(x) \,\, \mbox{solves
 (\ref{13})},
  \ee
  on compact subsets in $[0,L_0)$.
 Note that, for the desired result \ef{A1}, the following holds:
  \be
  \label{pp1}
   \tex{
    \frac 1{A^2(\t)} \sim \frac 1{\ln \t} \not \in L^1_s((1,\iy)),
     \quad \mbox{since} \quad \frac {{\mathrm d}s}{A^2(\t)} \sim
     {\mathrm d} \t,
    }
    \ee
    i.e., the perturbation is not integrable in $L^1$ that causes
    extra difficulties (though is a  typical feature of many  blow-up and extinction problems, \cite{AMGV}).
However,
   this stage of the analysis is
  rather standard, since the unperturbed equation
   \be
   \label{24N}
w_s= w^2(-w_{xxxx}+w)
 \ee
 is a gradient system (cf. \ef{10}), though fighting the
 perturbation terms in \ef{22N} can be rather technical (but
 hopefully not of a  principle issue). In standard blow-up theory
 such problems are tackled by  theorems on stability of
 $\o$-limit sets with respect to arbitrary small perturbations
 of the dynamical system, provided that the unperturbed problem is
 uniformly Lyapunov stable; see examples and  a large number of references
in \cite{AMGV}.

Anyway, this is not the case, where the main difficulty occur.
Namely, we need to match the behaviour \ef{22N}, \ef{23N} with
that corresponding to a proper one induced by the boundary
conditions at $x=L$ (or an FBP setting, which is hard to
distinguish from). Then we ought  to perform the linearization:
 \be
 \label{25N}
  \tex{
  w(x,s)=f(x)+Y(x,s)
  \LongA
  Y_s = \CC Y + 2 f Y(-Y^{(4)}+Y)+ \frac 1{2A^2} \, f+...\, ,
   }
   \ee
   where we omit all asymptotically smaller linear and nonlinear
   terms, which however can be important for matching.
   As usual, equation in \ef{25N} clearly shows that  full spectral theory
   of the linearized operator
 \be
 \label{op1}
 \CC=f^2(x)(-D^4_x+I) \inB x \in [0,L_0)
 \ee
    is necessary for extension of
   the solution  for $x \ge L_0$. Fortunately, it is symmetric in the weighted space
   $L^2_{\rho}$,
    with the weight $\rho=f^{-2}(x) \ge 0$, which is singular at
    the end point $x=L_0$.
 However, checking the crucial end point $x=L_0$ and using the
 asymptotics of the coefficient \ef{20}, yields Euler's equation
 in the spectral problem:
  \be
  \label{op2}
  C_1^2(L_0-x)^4 (-\psi^{(4)}+ \psi)+...= \l \psi.
   \ee
  Hence, linearly independent solutions  polynomial  with the characteristic equation:
  \be
  \label{op3}
  \psi(x) \sim (L_0-x)^m \LongA C_1^2 m(m-1)(m-2)(m-3)+\l=0.
   \ee
This shows that the necessary spectral theory, though covered by
classic theory of singular ordinary differential operators (see,
e.g., Naimark \cite{Nai1}), can be rather tricky.
 In particular, deficiency indices of the operator \ef{op1} are
 not easy to detect in $L^2_\rho$.
 More precisely,
depending on the boundary matching, this can require self-adjoint
extensions in both the cases of discrete and continuous spectra.
Infinitely oscillatory behaviour for complex values of $m= a+ \ii
b \in {\mathbb C}$ in \ef{op3} ($\l$ remains real) occurs
 close to $L_0^-$, where
 \be
 \label{op4}
 Y(x,s) \sim {\mathrm e}^{\l s} (L_0-x)^a \sin(b \ln(L_0-x)) \asA s \to +\iy;
  \ee
cf. the discussion of oscillatory issues in Section \ref{S1N}.
Thus, using proper spectra theory of \ef{op1} with eigenfunction
expansions of $w(x,s)$ of \ef{25N} over a discrete (better) or
continuous (most plausible for many problem settings) is the first
key concern here.

 Continuing  the spectral aspects of the expansion, we next assume
 that the spectrum of $\CC$ in \ef{op1} is continuous in a
 suitable setting and deficiency indices, so we may impose
 the following behaviour of the solution $Y$ of \ef{25N}:
  \be
  \label{op8}
   Y(x,s)= \mu(s) \Phi(x)+ ... \asA s \to +\iy.
    \ee
    Substituting this into \ef{25N} and
    naturally assuming that
     \be
     \label{a10}
     \tex{
     \mu'_s= \mu'_\t \t'_s \equiv - \frac {A'}{A^5}
      \LongA
    |\mu'(s)| \ll |\mu(s)| \forA s \gg 1
 }
    \ee
     (this is true  for the expectation \ef{A1}),
    we arrive at the following first rough balance:
     \be
     \label{op9}
      \tex{
      \mu'_s \Phi= \mu \CC \Phi + \frac 1{2A^2} f+... \LongA
      \mu(s)= \frac 1{2 A^2(\t)} \andA \CC \Phi+f=0.
    }
       \ee
 Again, a suitable solvability of the $\Phi$-equation in
 \ef{op9} requires posing proper boundary conditions at the
 singular end-point $x=L_0$ and satisfactory spectral theory of
 $\CC$.

 Meantime, we solve the ODE for $\Phi$ in \ef{op9}
 close to the singular end-point $x=L_0$:
  \be
  \label{eq1}
   \tex{
    - \Phi^{(4)}+ \Phi=- \frac 1{f(x)} = - \frac {1+o(1)}{C_1(L_0-x)^2}
 \asA x \to L_0^-.
 }
 \ee
 This yields the following non-regular asymptotic expansion:
  \be
  \label{eq2}
   \tex{
    \Phi(x)= - \frac 1{2C_1}\,z^2 \ln z+ A_1 + A_2 z+ A_3
    z^2 + A_4 z^3+... \asA z= L_0-x \to 0^+.
    }
    \ee
 Here we see the first clear appearance of the $\ln(L_0-x)$ in the
 framework of this standard matching and operator analysis:
 \be
 \label{eq3}
  \tex{
  \Phi(x) \sim  - \frac 1{2C_1}\,(L_0-x)^2 \ln (L_0-x)+A_1+A_2(L_0-x)+ ... \asA
  x \to L_0^-,
  }
  \ee
 which
 is indeed a reliable remnant of the asymptotic expansion \ef{19}
 for the log-TW profile $g(y)$ satisfying \ef{18}.
The integration constants $A_{1,2,3,4} \in \re$ in \ef{eq2} should
be  obtained by further matching (some of them can be arbitrary
and depend on initial data).


Next, in the R-II, where $w$ is uniformly small and hence $w^3$
can be neglected, the problem \ef{22N}, with the behaviour
\ef{op8} at the end-point $x=L_0$, takes the form
 \be
 \label{hh2}
  \left\{
  \begin{matrix}
  w_s= - w^2 w_{xxxx} + \mu(s) w+ ... \inB (L_0,L) \times \re_+, \qquad\qquad\,\, \ssk\\
 w =A_1 \mu(s), \,\, w_x=-A_2 \mu(s), \,\, x=L_0; \quad w=w_x=0, \,\, x=L.
  \end{matrix}
   \right.
   \ee
The precise behaviour of $w$ and the derivatives $w_x$ and
$w_{xx}$ at $x=L_0$ is defined by  \ef{op8}. Therefore, we add to
\ef{hh2} the following overdetermining condition
 \be
 \label{hh3}
  \tex{
  w_{xx}(x,s)= -\frac 1{C_1} \ln(x-L_0) (1+o(1)) \mu(s) \asA x \to
  L_0^+,
  }
  \ee
  which is expected to detect the unknown function $A(\t)$.
  Theoretically, once these have been achieved, extending
 the solution up to the original boundary point $x=L$ will define
 the necessary ODE-type equation for the unknown function $A(\t)$
 such as, very roughly,
 \be
 \label{hh1}
 A' \sim {\mathrm e}^{-A^2} \LongA A(\t) \sim \sqrt{\ln \t} \forA
 \t \gg 1.
  \ee

 The principle key difficulties start now, when we arrive at the problem
 \ef{hh2}. The point is that when $w(s,x)$ approaches 0, the
 problem becomes extremely sensitive and demand solving a number
 of linear and nonlinear eigenvalue problems, where we face in a
 full scale those unresolved local
 existence-uniqueness-entropy-etc. problems mentioned in Section
  \ref{S1N}.
Namely, even  the unperturbed problem
  \be
  \label{ww1}
 w_s=- w^2 w_{xxxx} \quad \mbox{on small solutions} \quad |w| \ll
 1
  \ee
is not well-understood at all. Note that the second-order
counterpart of \ef{par1},
 \be
 \label{ww3}
 w_s=w^2 w_{xx},
   \ee
    is also not easy, but by a  contact B\"acklund symmetry reduces to the
heat equation:
 \be
 \label{ww6}
  \tex{
 v_s=v_{yy} \whereA w(x,s)=\frac 1{\phi_x(x,s)} \andA
 \left\{
  \begin{matrix}y= \phi,\qquad\\ x=
v(y,s);
 \end{matrix}
 \right.
 }
  \ee
   see \cite[p.~78]{SGKM}
and references therein. Actually, such a homology
\ef{ww3}\,$\mapsto$\,\ef{ww6}  eventually creates a relation
\ef{hh1}, where the exponential term ${\mathrm e}^{-A^2}$ is
naturally associated with the Gaussian fundamental solution of the
operator $D_s-D^2_{y}$ in \ef{ww6}.
 Obviously, no such nonlocal symmetry  is available for the much
 more difficult PDE \ef{ww1}.

  Thus,
 the technicalities at this last stage are
 extremal and
 which deserves extra  special difficult analysis
that we
 cannot be performed in the present more formal paper (will require at least an extra
   dozen of sheets
 of hard calculations), so that, mathematically:
  $$
   \fbox{$
  \mbox{{\bf Main difficulty:} problems (\ref{hh2}) and (\ref{ww1})
  are not doable at the moment}
 $}
   $$
 (an open problem for the attentive Reader).
 Let us note that  spectral theory of
 generalized fourth-order non-symmetric Hermite operators
 (an operator pair)
 \cite{Eg4} can be key:
  \be
  \label{op10}
   \left\{
   \begin{matrix}
   \BB= -D_y^4 + \frac 14 \, y D_y + \frac 14 \, I, \ssk\\
\BB^*= -D_y^4 - \frac 14 \, y D_y \qquad\,\,
 \end{matrix}
  \right.
 \LongA
  \tex{
\s(\BB)=\s(\BB^*)=\big\{\l_k=- \frac k4, \, k \ge 0\big\}.
 }
 \ee

 As we have mentioned, exactly in this R-II,
 a delicate and specially deformed structure of log-TWs are crucial for
 creating the whole pattern.
Thus, the geometric approach of formal matching makes it possible
to avoid some (too) detailed refinement of the solution expansion,
and establishes  that the double-log factor can occur for a kind
of an ``envelope" of a wide class of blow-up solutions induced by
various functional settings.
 However, the behaviour is expected to be so difficult that we
 will not be surprised (e.g.,) if the resulting more correct
 calculation  provides us with a different double-log factor
   $
   \hat A(t)  \sim \big[ \ln|\ln(T-t)|\big]^{ 3/4}
   $
   instead of the standard one \ef{AA1}, though the topology
   matching  via \ef{23}, \ef{24} indicates that this is not
   the case.


Finally, let us point out again that  a full (formal) asymptotic
expansion technique also assumes the next matching of the log-TW
Region-II with the Boundary one close to $x=L$, where another (and
hopefully and typically much weaker than at $x=L_0$) boundary
layer of uniformly bounded solutions can occur. As usual, this is
expected to be easier (and less principal), but can be also
tricky, as everything concerning such nonlinear degenerate
higher-order operators. We recall that a proper and mathematically
justified matching of those manifolds of solutions via rigorous
compactness-like arguments seems not achievable for the
higher-order parabolic equations under consideration. Note again
that the analysis in \cite{AngVel1} for \ef{par1} was essentially
and inevitably based, among others, on Maximum and Comparison
Principle arguments, which are non-existent for higher-order
parabolic flows. Moreover, as we have pointed out in greater
detail in Section \ref{S1N}, the CP, IBVP, or any FBP settings for
\ef{p1} have not been well developed still, so we have to be very
careful in deriving asymptotic blow-up properties of some possibly
nonexistent and/or non-unique solutions.


\section{Briefly on  the sixth-order model: the log-log universality}
 \label{S2m3}

For the equation \ef{m3}, the ODE in \ef{p41} is
 \be
 \label{41}
  \tex{
   \frac 12 \, \th= \th^2(\th^{(6)}+\th),
   }
   \ee
so that nonexistence of a sufficiently smooth solution follows
from the identity (cf. \ef{p7})
 \be
 \label{42}
  \tex{
   \frac 12 \ln|\th|= \th^{(5)}\th'-\th^{(4)}\th''+ \frac
   12\,(\th'')^2+ \frac 12\, \th^2 +C.
   }
   \ee
The ODE for the log-TWs in Inner Region-II now takes the form
 \be
 \label{43}
  \tex{
  \frac 12\, g= g^2 g^{(6)}+... \LongA
   g(\eta) \sim \frac 1{2\sqrt 3}\, (-\eta)^3 \sqrt{\ln(-\eta)}
   \asA \eta \to -\iy.
   }
   \ee
Therefore, the same argument of matching with the stationary
structure in blow-up Inner Region-I governed by
 \be
 \label{44}
 f(x): \quad f^{(6)}+f=0 \inB (-L_0,L_0), \quad f=f'=f''=0 \,\,\,
 \mbox{at} \,\,\,x= \pm L_0
  \ee
  (hence with the cubic behaviour $f(x) \sim C_1(L_0-x)^3$ for $x \approx L_0$ that is
  well matched with the term $(-\eta)^3$ in \ef{43})
  yields as in \ef{24} the same log-log factor $A(t)$.

\section{No log-log for the PME--4 with source}
 \label{S3}

\subsection{Countable basic family of patterns by variational approach}

Consider the Cauchy problem for \ef{p3} in $\re \times (0,T)$.
Following the question \ef{Q1}, we then claim that there exist
weak compactly supported similarity solutions in the separable
variables:
 \be
 \label{25}
 \tex{
 u_{\rm S}(x,t)= \frac 1{\sqrt{T-t}}\, \th(x)
 \whereA \frac 12\, \th=-(\th^3)^{(4)}+\th^3 \inB \re.
 }
 \ee
 Changing for convenience the function,
  \be
  \label{26}
   \tex{
  \th(x)= \frac 1{\sqrt 2}\, F^{\frac 13}(x),
  }
  \ee
yields a more standard semilinear ODE with non-Lipschitz
nonlinearity,
 \be
 \label{27}
  F^{(4)}-F+F^{\frac 13}=0 \inB \re,
   \ee
which admit three obvious constant equilibria $0$ and $\pm 1$.
 Fortunately, this problem is variational, and the following
 result is proved by a combination of Lusternik--Schnirel'man
 category (genus) theory (L--S theory, for short) and Pohozaev's fibering approach of calculus
 of variations \cite{GMPSob}.

 \begin{proposition}
 \label{Pr.LS}
 The problem $\ef{27}$ admits at least a countable set of
 nontrivial compactly supported solutions.
  \end{proposition}

We present  a few comments concerning this principal result. We
look for critical points of the corresponding $C^1$-functional:
 \be
 \label{V1}
  \mbox{$
 {E}(F)= - \frac 12  \int (F'')^2 + \frac 12 \int
 F^2 -\frac 34 \, \int |F|^{\frac 43}.
  $}
  \ee
  In general, we have to look for critical points in $W^2_2(\ren)
  \cap L^2(\ren) \cap L^{4/3}(\ren)$.
  Bearing in mind 
 compactly supported solutions, we choose a sufficiently
  large  $R>0$  and consider the variational problem for (\ref{V1})
  in $W_{2,0}^2(B_R)$, $B_R=(-R,R)$, where we assume Dirichlet boundary
  conditions. It is next proved that any  solutions satisfying
  $ F(y) \to 0$ as $|y| \to \infty$
 is compactly supported in $\re$.

 Thus, the functional (\ref{V1}) is $C^1$ and is
uniformly differentiable and weakly continuous, so we can apply
classic
 Lusternik--Schnirel'man (L--S) theory of calculus of variations
\cite[\S~57]{KrasZ} in the form of the fibering method \cite{Poh0,
PohFM}.
 Namely, following  L--S theory and the fibering
approach \cite{PohFM}, the number of critical points of the
functional (\ref{V1}) depends on the {\em category} (or {\em
genus}) of the functional subset on which fibering is taking
place. Critical points of ${E(F)}$ are obtained by
 {\em spherical
fibering}
 \be
 \label{f1}
 F= r(v) v \quad (r \ge 0),
  \ee
  where $r(v)$ is a scalar functional, and $v$ belongs to a subset
  in  $W_{2,0}^2(B_R)$ given as follows:
   \be
   \label{f2}
    \mbox{$
    {\mathcal H}_0=\bigl\{v \in W_{2,0}^2(B_R): \,\,\,H_0(v)
     \equiv  -  \int (v'')^2 +  \int
 v^2 =1\bigr\}.
    $}
    \ee
The new functional
 \be
 \label{f3}
  \mbox{$
H(r,v)= \frac 12 \, r^2 - \frac 34\, r^{\frac 43} \int |v|^{\frac
43}
 $}
  \ee
 has the absolute minimum point, where
 \be
 \label{f31}
  \mbox{$
 H'_r \equiv r-  r^{\frac 13} \int |v|^{\frac 43} =0
  \,\,\Longrightarrow \,\,
   r_0(v)=\bigl(\int |v|^{\frac 43}\bigr)^{\frac 32}.
   $}
   \ee 
   We then obtain the following functional:
 \be
 \label{f400}
 \mbox{$
  \tilde H(v)= H(r_0(v),v)=- \frac 14 \,  r_0^2(v)
  \equiv - \frac 14 \, \big( \int |v|^{\frac 43}\big)^{3}.
   $}
    \ee
 The critical points of the functional (\ref{f400}) on the
 set (\ref{f2}) coincide with those for
 \be
  \label{f4}
   \mbox{$
    \tilde H(v)= \int |v|^{\frac 43},
    $}
    \ee
 so
  we arrive at even, non-negative, convex,
 and uniformly differentiable functional, to which
L--S  theory applies, \cite[\S~57]{KrasZ}; see also
\cite[p.~353]{Deim}.
   Following \cite{PohFM}, searching
  for critical points of $\tilde H$ on the set ${\mathcal H}_0$
  one needs to estimate the category $\rho$
  of the set ${\mathcal H}_0$.
The details on this notation and basic results can be  found in
Berger \cite[p.~378]{Berger}.

It follows that, by this variational construction, $F$ is an
eigenfunction satisfying
 $
 - F^{(4)} + F - \mu F^{\frac 13}=0,
  $
where $\mu >0$ is Lagrange's multiplier. Then scaling $F \mapsto
\mu^{3/2} F$ yields the original equation in (\ref{27}).

\ssk

For further discussion of geometric shapes of patterns, it is
convenient to recall that utilizing Berger's version
\cite[p.~368]{Berger} of this minimax analysis of L--S category
theory \cite[p.~387]{KrasZ},  the critical values $\{c_k\}$ and
the corresponding critical points $\{v_k\}$ are given by
 \be
 \label{ck1}
  \mbox{$
 c_k = \inf_{{\mathcal F} \in {\mathcal M}_k} \,\, \sup_{v \in {\mathcal
 F}} \,\, \tilde H(v),
  $}
  \ee
where  ${\mathcal F} \subset {\mathcal H}_0$ are  closed sets,
 and
 ${\mathcal M}_k$ denotes the set of all subsets of the form
  $
  B S^{k-1}
\subset {\mathcal H}_0,
 $
 where $S^{k-1}$ is a suitable sufficiently
smooth $(k-1)$-dimensional manifold (say, sphere) in ${\mathcal
H}_0$ and $B$ is an odd continuous map.
 Then each member of ${\mathcal M}_k$ is of  genus at least $k$
 (available in ${\mathcal H}_0$).
   It is also important to remind that the
definition of genus \cite[p.~385]{KrasZ} assumes  that
$\rho({\mathcal F})=1$, if no {\em component} of ${\mathcal F}
\cup {\mathcal F}^*$, where
 $
 {\mathcal F}^*=\{v: \,\, v^*=-v \in {\mathcal F}\},
 $
 is the {\em reflection} of ${\mathcal F}$ relative to 0,
 contains a pair of antipodal points $v$ and $v^*=-v$.
 Furthermore, $\rho({\mathcal F})=n$ if each compact subset of
${\mathcal F}$ can be covered by, minimum, $n$ sets of genus one.

According to (\ref{ck1}),
 $
 c_1 \le c_2 \le ... \le c_{l_0},
 $
 where $l_0=l_0(R)$ is the category of ${\mathcal H}_0$ (see an estimate
 below) such that
  \be
  \label{l01}
  l_0(R) \to + \infty \quad \mbox{as} \quad R \to \infty.
  \ee
  Roughly speaking,
since the dimension of the sets ${\mathcal F}$ involved in the
construction of ${\mathcal M}_k$ increases with $k$, this
guarantees that the critical points delivering critical values
(\ref{ck1}) are all different.


  It follows from (\ref{f2}) that the category
$l_0=\rho({\mathcal H}_0)$  of the set ${\mathcal H}_0$ is equal
to the number (with multiplicities) of the eigenvalues $\l_k>-1$,
  \be
  \label{pp1nn}
 \mbox{$
  l_0= \rho({\mathcal H}_0)= \sharp \{ \l_k > -1\}
 $}
  \ee
 of the linear bi-harmonic operator $- D^4<0$,
 \be
 \label{f55}
 -  \psi_k^{(4)}= \l_k \psi_k, \quad \psi_k \in W^2_{4,0}(B_R);
  \ee
  see \cite[p.~368]{Berger}.
  Since the dependence of the spectrum on $R$ is, obviously,
   \be
   \label{f56}
   \l_k(R)= R^{-4} \l_k(1), \quad k=0,1,2,... \, ,
    \ee
we have that the category $\rho({\mathcal H}_0)$ can be
arbitrarily large for $R \gg 1$, and (\ref{l01}) holds.


\subsection{On total variety of patterns: numerical evidence}

Actually, the total variety of possible solutions is not exhausted
by the categories of L--S theory, which we will explain below by
presenting clear numerical evidence.

\begin{figure}
\centering
\includegraphics[scale=0.75]{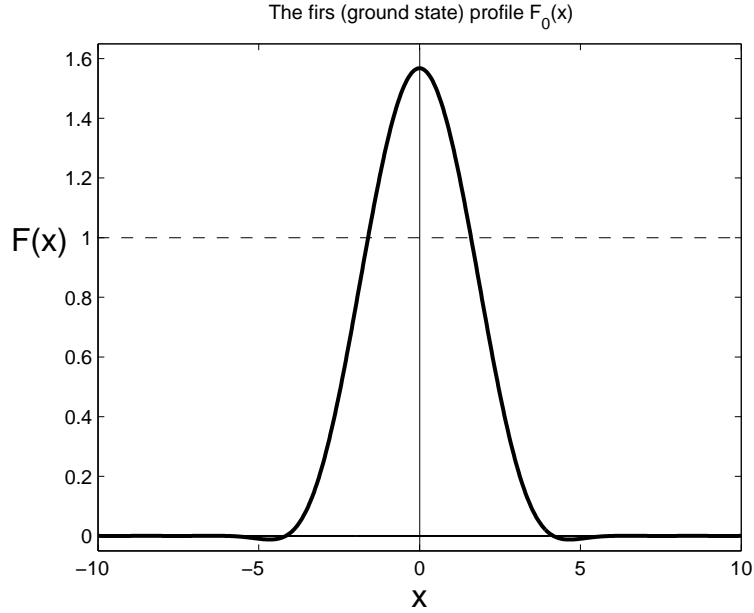} 
\vskip -.3cm \caption{\small The first, ground state solution of
the ODE \ef{27}.}
   \vskip -.3cm
 \label{F2}
\end{figure}

 Figure \ref{F2} demonstrates the first, and actually, the {\em ground
 state} profile $F_0(x)$. It has a typical shape of a ground
 state, as a critical point of the functional delivering the
 absolute extremum. Note that unlike the classic second-order case
 of the ground state for \cite{Coff72} (this problem is key for
 the critical NLSE)
  \be
  \label{gr1}
  \D u - u + u^3 =0 \inB \ren,
   \ee
   which is {\em strictly positive} (with exponential decay at infinity) and
   is unique up to translations, the ground state $F_0$ for
   \ef{27} is oscillatory and of changing sign near {\em finite
   interfaces}; this is seen from Figure \ref{F2}.  The oscillatory
   structure of solutions is rather involved and is described in
   \cite{GMPSob}; 
    see also similar details in \cite{Gl4}.

\begin{figure}
\centering
\includegraphics[scale=0.75]{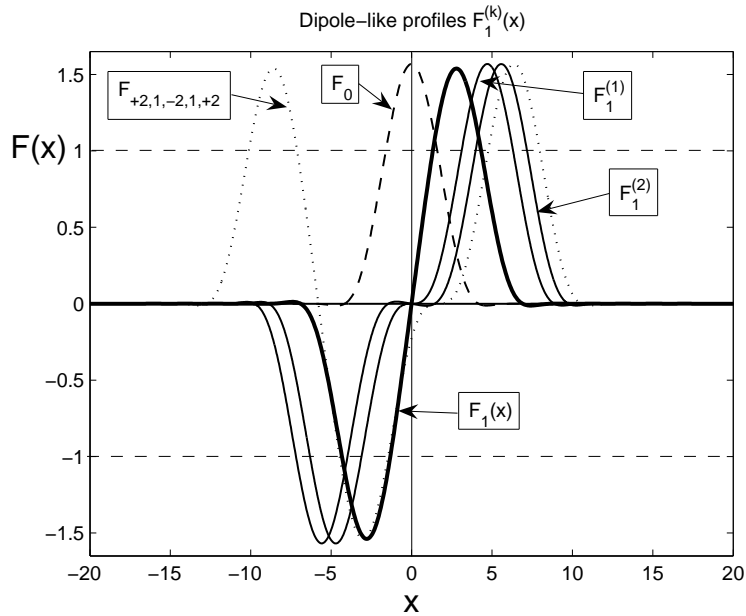} 
\vskip -.3cm \caption{\small Some of the dipole-like solutions
$F_1^{(k)}(x)$, $k=0,1,2$ of the ODE \ef{27}.}
   \vskip -.3cm
 \label{F3}
\end{figure}

In Figure \ref{F3}, we show the second dipole-like profiles, where
$F_1(x)$ denoted by the boldface line is the basic one
corresponding to the L--S critical value. In addition, there
exists a countable (this is not proved still) family of dipoles
$\{F_1^{(k)}\}$, which differ by the structure of the internal
zero close to the origin (in general, $k$ stands for the number of
transversal zeros there, but this is not enough to uniquely
identifying the pattern). By the dotted line, we denote another
profile from the next family $\{F_2^{(k)}\}$; see below.

\begin{figure}
\centering
\includegraphics[scale=0.75]{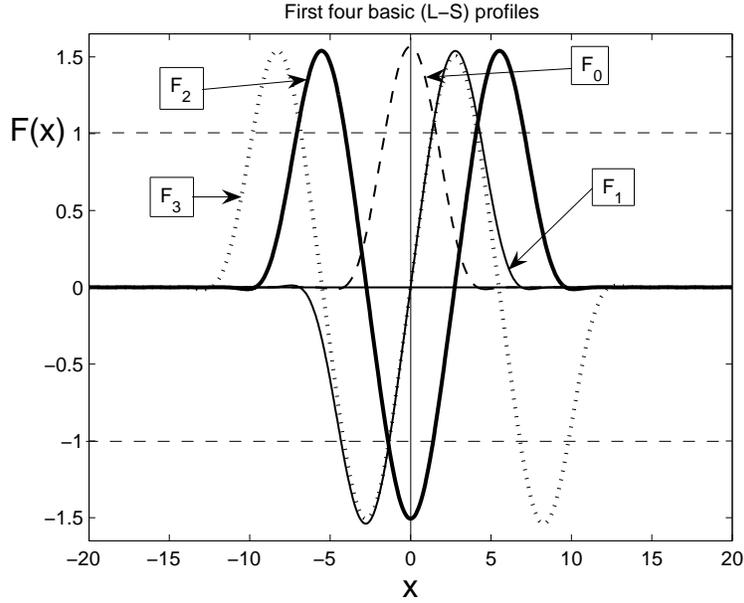} 
\vskip -.3cm \caption{\small Four basic L--S patterns of the ODE
\ef{27}.}
   \vskip -.3cm
 \label{F3New}
\end{figure}

Figure \ref{F3New} continues explaining further basic L--S
patterns, where we show $F_0$, $F_1$, $F_2$, and $F_3$ (the dotted
line). It is clearly seen that each $F_k$ has precisely $k$
``dominant" transversal zeros inside the support, which well
corresponds to Sturm-like principle (not applicable here in the
rigorous sense, since all the solutions are oscillatory and have
infinitely many sign changes near interfaces).

\begin{figure}
\centering
\includegraphics[scale=0.75]{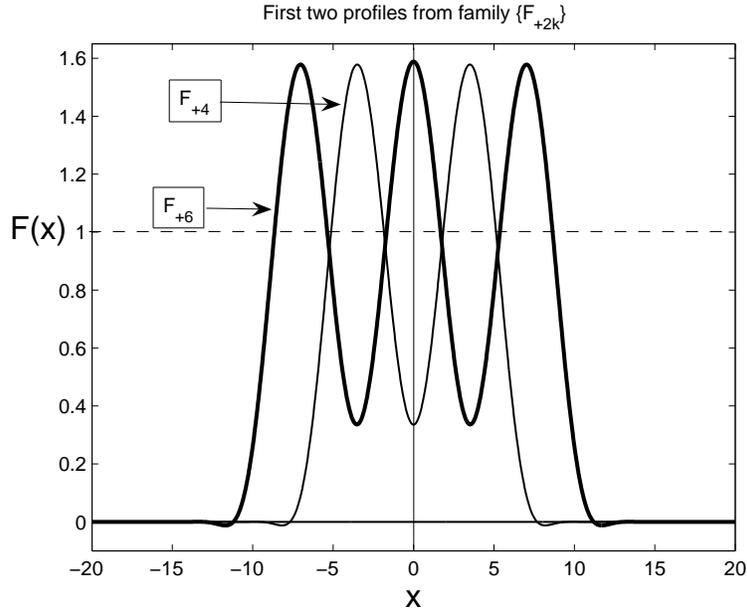} 
\vskip -.3cm \caption{\small  Solutions $F_{+4}(x)$ and
$F_{+6}(s)$ of the ODE \ef{27}.}
   \vskip -.3cm
 \label{F4}
\end{figure}

\begin{figure}
\centering
\includegraphics[scale=0.75]{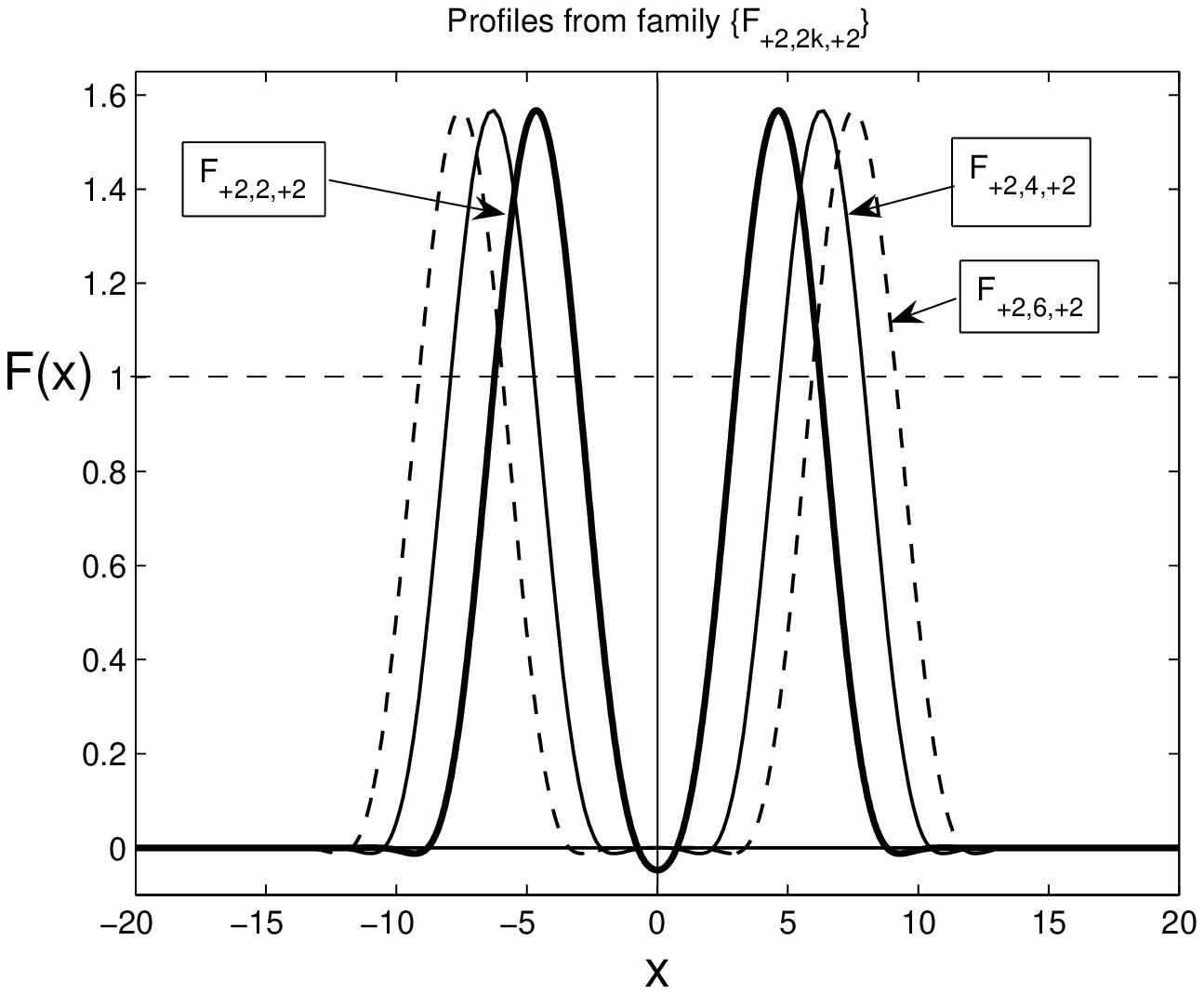} 
\vskip -.3cm \caption{\small  Solutions $\{F_{+2,2k,+2}(x)\}$  of
the ODE \ef{27}.}
   \vskip -.3cm
 \label{F5}
\end{figure}

In Figure \ref{F4}, non L-S profiles from the family $\{F_{+2k}\}$
are presented. Figure \ref{F5} shows some profiles from the family
$\{F_{+2,2k,+2}(x)\}$, which are also not expected to correspond
to any L--S category/critical value.

\begin{figure}
\centering
\includegraphics[scale=0.75]{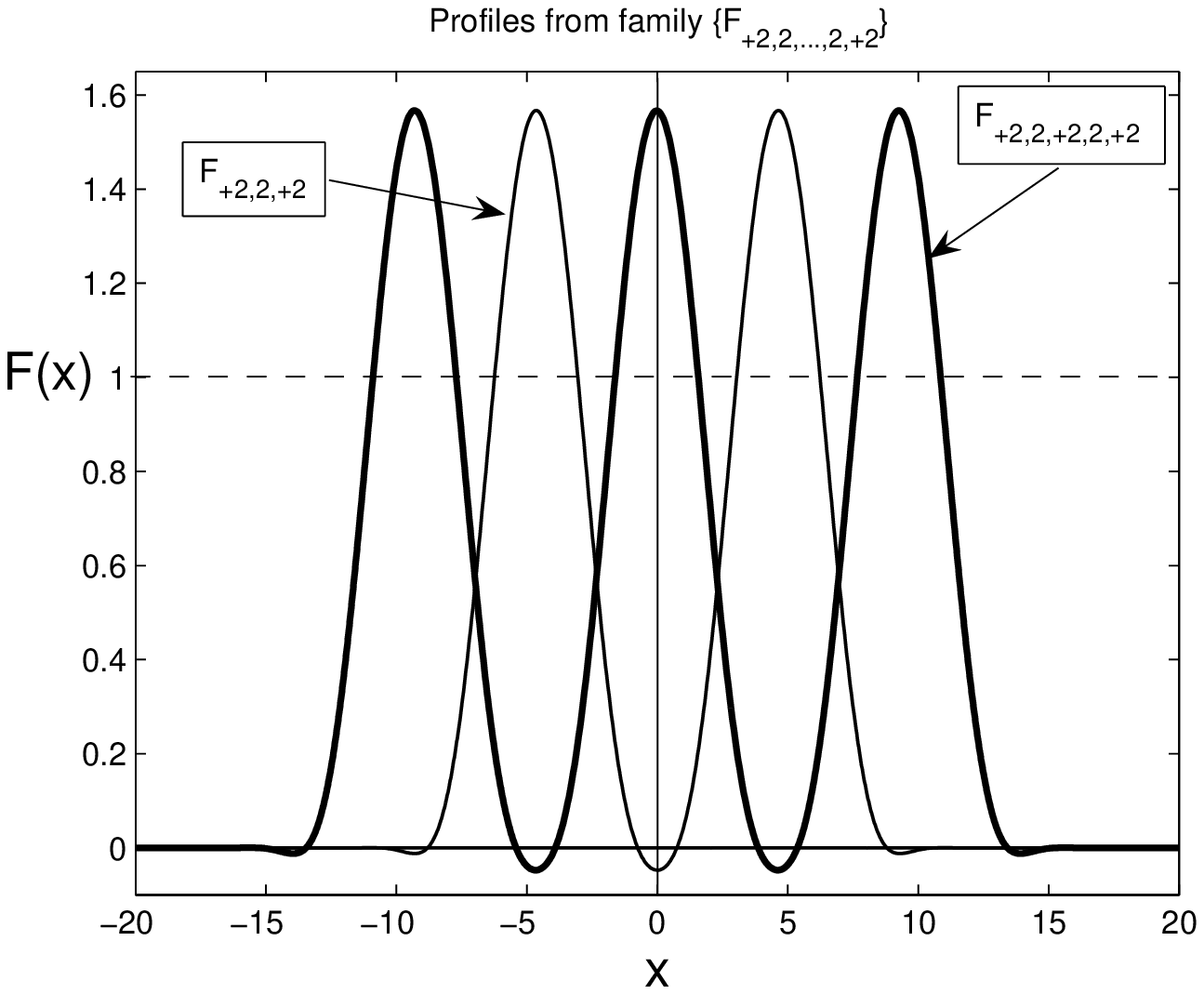} 
\vskip -.3cm \caption{\small  Solutions $\{F_{+2,2,...,2,+2}(x)\}$
of the ODE \ef{27}.}
   \vskip -.3cm
 \label{F6}
\end{figure}

First profiles from the family $\{F_{+2,2,...,2,+2}(x)\}$ (also
non L--S) are shown in Figure \ref{F6}.

\begin{figure}
\centering
\includegraphics[scale=0.75]{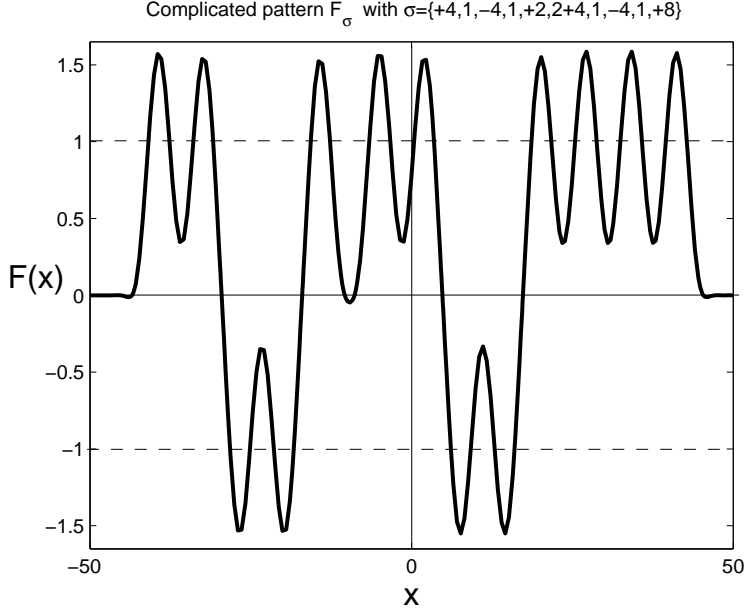} 
\vskip -.3cm \caption{\small  Solution
$F_{+4,1,-4,1,+2,2,+4,1,-4,1,+8}(x)$ of the ODE \ef{27}.}
   \vskip -.3cm
 \label{F7}
\end{figure}

Finally, following \cite{GMPSob}
and using the above rather
simple families of patterns, we claim that a pattern (possibly, a
class of patterns) with an arbitrary multiindex of any length
 \be
 \label{mm1}
 \s=\{\pm \s_1, \s_2, \pm \s_3, \s_4,..., \pm \s_l\}
  \ee
  can be constructed. Here, as above, each $\pm \s_k$
  stands for the total number of successive  intersections with the current
 equilibrium $\pm 1$, while $\s_k$ counts that with the trivial
 equilibrium $0$.
E.g.,
  in Figure \ref{F7}, we show
 a single complicated profile $F_\s(x)$
  with the index
   \be
   \label{ch1}
   \s=\{+4,1,-4,1,+2,2,+4,1,-4,1,+8\}.
    \ee

Actually,  the multiindex (\ref{mm1}) can be rather arbitrary
(with some natural restrictions on even local intersection
numbers) and then takes  finite parts of any admissible
non-periodic fraction. Overall, this means {\em chaotic features}
 of the whole family of solutions $\{F_\s\}$.
 These chaotic types of behaviour are known for other fourth-order ODEs with
coercive smooth operators, \cite[p.~198]{PelTroy}.

\section{No log-log for the TFE--4 with source}
 \label{S4}

For the TFE--4 \ef{p4}, the similarity substitution as in \ef{25}
yields another ODE:
 \be
 \label{k1}
  \tex{
   \frac 12\, \th=-(\th^2 \th''')'+ \th^3 \inB \re.
   }
   \ee
To get equilibria $\pm 1$, we perform the change:
 \be
 \label{k2}
 \tex{
 \th(x)= \frac 1{\sqrt{2}} \, F(x) \LongA (F^2 F''')'=F^3-F \inB
 \re.
 }
 \ee

\begin{figure}
\centering
\includegraphics[scale=0.75]{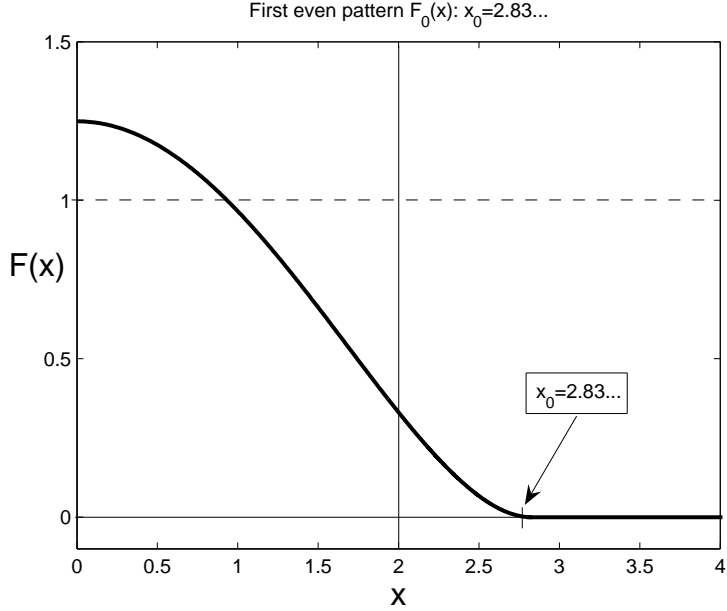} 
\vskip -.3cm \caption{\small  The first basic nonnegative solution
$F_0(x)$ of the ODE in \ef{k2}.}
   \vskip -.3cm
 \label{F8}
\end{figure}

   Unlike the one in \ef{25} and \ef{27}, the ODE in \ef{k2} does not admit a
 variational formulation. Nevertheless, there exist rather
 standard shooting arguments for detecting necessary similarity
 blow-up profiles; see references in \cite[\S~3]{Gl4}. In
 particular, the following asymptotic behaviour is known close to
 finite interfaces\footnote{These are the {\em maximal regularity}
 solutions admitted by the ODE \ef{k2}; see \cite{Gl4}.}
  \be
  \label{k3}
   \tex{
   F(x)= \frac 12\,(x_0-x)_+^2 \big[\sqrt{|\ln(x_0-x)|}+C+...\big] \asA x
   \to x_0,
   }
   \ee
 where $C \in \re$ is an arbitrary constant. Overall, the asymptotic bundle
 \ef{k3} comprises {\em two} parameters $\{x_0,C\}$, which are
 expected to be sufficient to shoot also {\em two} conditions at
 the origin:
  \be
  \label{k4}
  F'(0)=F'''(0)=0 \,\,\,(\mbox{symmetry}) \quad \mbox{or}
  \quad  F(0)=F''(0)=0 \,\,\,(\mbox{anti-symmetry}).
   \ee
 Figure \ref{F8} shows the first nonnegative even profile
$F_0(x)$ for $x>0$ with the symmetry conditions in \ef{k4} and
 the interface
 $
 x_0=2.83...\, .
  $
 Hence, the answer to \ef{Q1} is log-log $\not \exists$.

\section{Quasilinear wave and nonlinear dispersion equations}
 \label{S5}

\subsection{The QWE--4}

For \ef{QWE1}, the similarity separable solution \ef{p41} is
slightly different, but the ODE remains the same:
 \be
 \label{q41}
 \tex{
 u_{\rm S}(x,t)= \frac 1{T-t}\, \th(x)
  \whereA
 2 \th= \th^2(-\th^{(4)} +  \th).
  }
  \ee
Hence,  Proposition \ref{Pr.Non} holds, so we need the log-TWs for
matching:
 \be
 \label{q6}
  \tex{
  u_{\rm log}(x,t)= \frac 1{T-t}\, g(\eta), \,\, \eta=x+\l
  \ln(T-t)\, \Longrightarrow \, 2g-3\l g'+\l^2 g''= g^2(-g^{(4)}+g).
   }
   \ee
However, the necessary asymptotics
 is guaranteed by
keeping just two terms:
 \be
 \label{q7}
 2g=-g^2 g^{(4)}+...\, ,
  \ee
  so that the only correction to \ef{19} is the multiplier $\sqrt
  2$.
The rest of the analysis remains the same as at the end of Section
\ref{S2}.

 For the divergent wave model
 \be
 \label{dv1}
 u_{tt}=-(u^3)_{xxxx} + u^3
  \ee
the log-log factor is nonexistent, since it admits separable
compactly supported solutions
 \be
 \label{dv2}
  \tex{
  u_{\rm S}(x,t)= \frac 1{T-t}\, \th(x) \LongA 2 \th=
  -(\th^3)^{(4)}+\th^3,
  }
  \ee
  which on scaling reduces to that in \ef{25}.

\subsection{The NDE--4: non-divergent model}


For \ef{NDE1}, the similarity separable solution \ef{p41} is
slightly different, but the ODE remains the same:
 \be
 \label{q41n}
 \tex{
 u_{\rm S}(x,t)= (T-t)^{-\frac 13}\, \th(x)
  \whereA
  \frac {1} 3 \,  \th= \th^{3}(\th''' +  \th).
  }
  \ee
  Instead of \ef{p7}, the nonexistence conclusion is governed by
  the following ``local monotonicity" identity:
 \be
 \label{po1n}
  \tex{
  - \frac {1} 3 \, \frac 1 \th=  \th' \th''- \int (\th'')^2 + \frac
  12\, \th^2,
   }
   \ee
   so that sufficiently  smooth  solutions with $\th=0$ at $\pm L$ (and $\th \in H^2$)
    do not exist.

As the analogy to \ef{14}, consider the stationary equation:
 \be
 \label{st1n}
 f(x): \quad f'''+f=0 \quad \mbox{on} \quad (-L_0,L_0), \quad
 f(-L_0)=f'(-L_0)=0, \quad f(L_0)=0,
  \ee
  where the zero at $x=L_0$ is assumed to be transversal with the
  behaviour (cf. the smoother one \ef{20})
   \be
   \label{st2n}
   f(x)=C_1(x-L_0)(1+o(1)) \asA x \to L_0^- \quad (C_1>0).
    \ee

According to the problem \ef{st1n}, to avoid extra difficulties
with rather obscure consequences, we consider \ef{NDE1} in
$(-L_0,L) \times (0,T)$, with $L>L_0$ (possibly, $L \gg L_0$ to
avoid existence of other stationary profiles) and the same
boundary conditions at the left-hand end point,
 \be
 \label{st3n}
 u=u_x=0 \quad \mbox{at} \quad x=-L_0 \andA u=0 \quad \mbox{at} \quad
 x=L.
  \ee
 Therefore, the log-log perturbation will penetrate into Inner
 Region-I from the right-hand singular point $x=L_0$, where the
 stationary solution \ef{st1n} vanishes thus creating an internal
 singular layer to be resolved by using slowly moving log-TWs.

  Thus,  the log-TW ansatz yields the ODE
 \be
 \label{q61n}
  \tex{
  u_{\rm log}(x,t)= (T-t)^{-\frac 13} g(\eta), \,\, \eta=x+\l
  \ln(T-t) \Longrightarrow  \frac 13\, g- \l g'= g^3(g'''+g).
   }
   \ee
 Close to the necessary point $x=L_0$, we
keep as usual  two terms that yields the desired behaviour:
 \be
 \label{q7n}
  \tex{
  \frac 13\, g=g^3 g'''+... \LongA g(\eta) \sim
  (-\eta)[\ln(-\eta)]^{\frac 13} \asA \eta \to -\iy.
  }
  \ee
Similar to \ef{23}--\ef{24}, matching the linear structure
$(-\eta)$ in \ef{q7n} with the linear one $\sim(L_0-x)$ in
\ef{st2n} yields
 \be
 \label{q8n}
  A(t) \sim \big[\ln|\ln(T-t)|\big]^{\frac 13}(1+o(1)) \asA t \to T^-.
   \ee
Observe changing $\sqrt{(\cdot)}$ into the cubic root due to cubic
nonlinearity $u^3$ in \ef{NDE1} instead of $u^2$ in other models.
The choice of the quartic equation \ef{NDE1} was in fact generated
by the possibility of such a transversal matching of ``linear"
structures. Some aspects of such a matching are still unclear and
deserve further study.

\subsection{No double-log in a divergent NDE}

\begin{figure}
\centering
\includegraphics[scale=0.65]{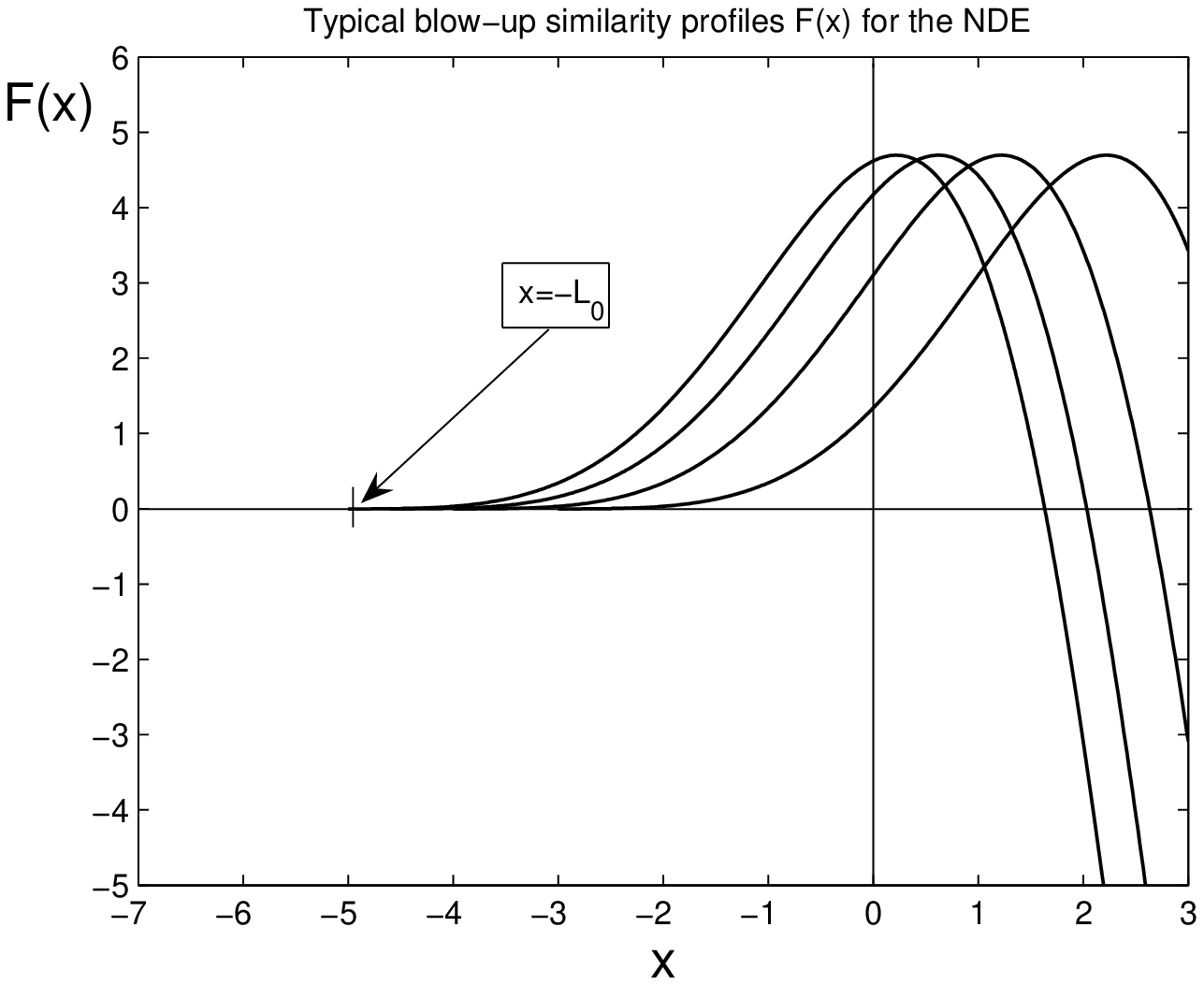} 
\vskip -.3cm \caption{\small  Blow-up patterns $F(x)$ satisfying
\ef{re1} with the behaviour \ef{re2}.}
   \vskip -.3cm
 \label{FLast}
\end{figure}

 For the corresponding  divergent NDE,
 \be
 \label{dv1n}
 u_{t}=(u^4)_{xxx} + u^4,
  \ee
 answering \ef{Q1}, one can construct  separable blow-up solutions
 \be
 \label{dv2n}
  \tex{
  u_{\rm S}(x,t)= (T-t)^{-\frac 13} \,  \th(x) \LongA  \frac 13\, \th=
  (\th^4)'''+\th^4.
  }
  \ee
Looking for nonnegative solutions $\th \ge 0$, on rescaling,
yields
 \be
 \label{re1}
 \th= 3^{-\frac 34} F \LongA F'''+F- F^{\frac 14}=0 \quad
 \mbox{in} \quad (-L_0,L_0).
  \ee
Unlike the ODE \ef{27}, this one admits solutions with a single
interface (free-boundary) point $x=-L_0$ with the behaviour
 \be
 \label{re2}
 F(x) =(24)^{- \frac 43}(L_0+x)_+^4(1+o(1)) \asA x \to -L_0^+.
  \ee
  In other words, we can set $F(x) \equiv 0$ for all $x <- L_0$.
  Therefore, moving such a profile $F(x)$, one can always satisfy
  the boundary condition at $x=L$ in \ef{st1n} for any $L>0$.
  Typical solutions of the ODE \ef{re1} are shown in Figure
  \ref{FLast}, where positive humps can serve as blow-up patterns.
   Note that all such solutions are not oscillatory
  near finite interfaces, unlike those for the TFEs \cite{Gl4}, where
  another third-order {\em oscillatory} ODE occurs.
 Thus, the similarity law \ef{dv2n} describes blow-up in the
 divergent model \ef{dv1n}.

\section{Final conclusion: blow-up log-log is universal
in  PDE theory}

The above study allows us to
 fix the following conclusion: rather surprisingly,
    \be
    \label{SSS}
   \fbox{$
    \mbox{
   the factor $\sqrt{\ln|\ln(T-t)|}$ has a clear {\bf universality} in
   blow-up
   }
   $}
   \ee
  for different classes of higher-order (non-divergent) nonlinear evolution PDEs.
  The
 $\ln\,\ln$-factor also occurs in boundary regularity (Petrovskii-type)  analysis
 for equations  \cite{GalPet2m}
  $$
  u_t=-u_{xxxx}, \quad u_t=u_{xxx}, \quad u_t=-(u^3)_{xxxx}, \quad \mbox{etc.}
  $$
 It would be  important, on the basis on the above discussion
 of various linear and nonlinear PDEs and by adding new asymptotic
 phenomena of necessity, to explain, in a more unified way, and to derive, formally or more justified,
such a common
  matched
 asymptotic $\sqrt{\log\,\log}$-criterion.



\end{document}